\documentclass{svjour3}
\usepackage[utf8]{inputenc}
\usepackage{bbding}
\usepackage{amsmath,amssymb}
\usepackage[dvipsnames]{xcolor}
\usepackage{graphicx}
\usepackage{cite}
\usepackage{subfig}
\usepackage{placeins}
\usepackage{tikz}
\usepackage{algorithm}
\usepackage{algpseudocode}
\usetikzlibrary{decorations.pathmorphing,patterns}

\let\div\undefined\DeclareMathOperator{\div}{div}

\newcommand{\blue}[2][blue]{\emph{\textcolor{#1}{#2}}}

\tikzset{
	pline/.style={
		every plot/.style={
			mark=x,
			mark options={
				black,
				thick,
			},
			mark size=4pt,
		},
		very thick,
		black,
	},
}

\title{Ensemble Kalman Inversion method for an inverse problem in soil-structure interaction}
\author{L. Scandurra$^1$}

\institute{{\Envelope} L. Scandurra \at
			\email{leonardo.scandurra@units.it} \at \at
            $^{1}$Università degli Studi di Trieste, Department of Mathematics and Geosciences, Trieste, Italy
}

\date{Received: date / Accepted: date}

\begin{document}
\titlerunning{EKI method for inverse problems}
\maketitle

\begin{abstract}
The interaction between the foundation structures and the soil has been developed for many engineering applications. For the determination of the stress in foundation structure it is needed to determine the influence of the stiffness of soil with respect to the displacement w of the deformable plate (direct problem), and viceversa, how the stiffness of the foundation structure affects the resulting subsidence (inverse problem). In this paper, we deal with the Winkler mathematical model and propose to use an efficient Ensemble Kalman Inversion scheme (EKI) that regularizes iteratively the ill-posedness of the inverse problem. It is a regularizing optimizer used in Bayesian inverse problems that samples particles in pseudo-time introducing a motion due to the movement of these particles. The EKI algorithm converges to the solution of an optimization problem that minimizes the objective function. In this context we show how to reconstruct the Winkler subgrade reaction coefficient of a rectangular thin plate loaded with an existing building by using the EKI methodology combined by the finite difference method (FDM) to discretize the biharmonic operator of the governing equations.

\end{abstract}

\section{Introduction}\label{intro}
The soil-structure interaction (SSI) is a very widespread problem in the geotechnical and structural engineering fields. In the continuum elasticity approach for the analysis of a plate on elastic foundation, the Winkler model (1867) \cite{winkler1867lehre} has been adopted to determine the subgrade reaction coefficient $k$, $k>0$ that describes a linear relationship between the reactive pressure $p$ at any arbitrary point $x$ of the plate $\Omega$ and the deflection $w$ of the underlying soil:
\begin{equation}
	p(x) = kw(x),\quad x\in\Omega.
	\label{eq:reactive_pressure}
\end{equation}
The mechanical behavior of the soil is highly complex, due to its non-linearity, shear and strain stresses, anisotropic and heterogeneous nature, so that in the most of cases, the subgrade reaction coefficient is replaced by using a simplified model. In the Winkler model the foundation is assumed to be elastic and laying on an independent linear fictitious springs system whose stiffness corresponds to the subgrade reaction coefficient $k$. One of the main difficulties encountered for reconstructing $k$ is due to existing buildings that make measurements not directly accessible for experiments. This work refers to the recent papers \cite{alessandrini2015global} where an inverse problem is considered in order to determine the Winkler coefficient $k$. We briefly describe the mechanical model:

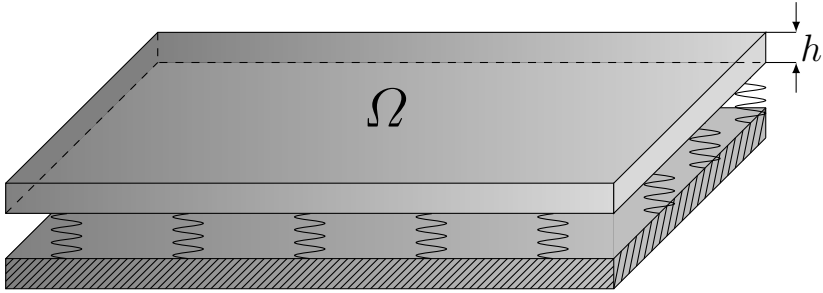
\begin{figure}[!ht]
	\centering
	\begin{tikzpicture}[scale=2]
			\fill[left color=gray, right color=lightgray] (-1,-1.3) -- (3,-1.3) -- (3,-1) -- (-0.7,-1) -- cycle;
			\fill[lightgray] (3,-1.3) -- (4,-0.3) -- (3.5,-0.3) -- (3,-1) -- cycle;
			\fill[lightgray] (3,-1.5) -- (4,-0.5) -- (4,-0.3) -- (3,-1.3) -- cycle;
			\fill[left color=gray, right color=lightgray] (-1,-1.5) rectangle (3,-1.3);
			\draw [decoration={aspect=0, segment length=1.7mm, amplitude=2mm,coil},decorate](3.3,-1) -- (3.3,-0.7);
			\draw [decoration={aspect=0, segment length=1.7mm, amplitude=2mm,coil},decorate](3.6,-0.7) -- (3.6,-0.4);
			\draw [decoration={aspect=0, segment length=1.7mm, amplitude=2mm,coil},decorate](3.9,-0.4) -- (3.9,-0.1);
			\fill[left color=gray, right color=gray!30] (0,0) rectangle (4,0.2);
			\fill[gray] (-1,-1) -- (0,0) -- (0,0.2) -- (-1,-0.8) -- cycle;
			\fill[left color=gray, right color=gray!30] (-1,-1) -- (3,-1) -- (4,0) -- (0,0) -- cycle;
			\draw[left color=gray, right color=gray!30] (-1,-1) rectangle (3,-0.8);
			\draw[fill=gray!30] (3,-1) -- (4,0) -- (4,0.2) -- (3,-0.8) -- cycle;
			\draw (4,0.2) -- (0,0.2) -- (-1,-0.8);
			\draw[dashed] (0,0.2) -- (0,0) -- (4,0);
			\draw[dashed] (-1,-1) -- (0,0);
			\draw (4,0) -- (4.2,0);
			\draw (4,0.2) -- (4.2,0.2);
			\draw[-latex] (4.2,-0.2) -- (4.2,0);
			\draw[latex-] (4.2,0.2) -- (4.2,0.4);
			\node at (4.3,0.1) {\Large{$h$}};
			\node at (1.5,-0.3) {\huge{$\Omega$}};
			\foreach \x in {-0.6,0.2,1,1.8,2.6}
			{  
					\draw [decoration={aspect=0, segment length=1.7mm, amplitude=2mm,coil},decorate](\x,-1.3) -- (\x,-1);
				}
			\draw (3,-1.3) -- (4,-0.3);
			\draw (-1,-1.5) rectangle (3,-1.3);
			
			\draw (3,-1.5) -- (4,-0.5);
			\draw (4,-0.5) -- (4,-0.3);
			
			\draw (-1,-1.3) -- (-0.7,-1);
			\draw (4,-0.3) -- (3.7,-0.3);
			\foreach \x in {0,0.05,...,3.85}
			{ 
					\draw (-1+\x,-1.5) -- (-0.8+\x,-1.3);
				}
			\foreach \x in {0,0.05,...,0.15}
			{
					\draw (-1,-1.35-\x) -- (-0.95+\x,-1.3);
					\draw (3,-1.35-\x) -- (2.85+\x,-1.5);
				}
			\foreach \x in {0,0.05,...,1}
			{
					\draw (3+\x,-1.5+\x) -- (3.05+\x,-1.25+\x);
				}
		\end{tikzpicture}
	\caption{A clamped rectangular thin plate resting on an elastic foundation.}
	\label{fig:domain}
\end{figure}
\noindent The domain $\Omega\times[-h/2,h/2]$ (see Fig. \ref{fig:domain}) is a clamped rectangular thin elastic plate $\Omega$ with uniform thickness $h$ subjected to a concentrated force $f\delta(P_0)$, with $f\in\mathbb{R}$, $f>0$ and $\delta(P_0)$ defined as a Dirac delta applied at a specific internal point $P_0$ of the plate $\Omega$. We are able to express the stresses on the plate as a function of a deflection $w$ of the plate, because 
the point on which the load is applied perpendicularly before bending, keeps its perpendicularity after bending and the normal stresses in the transverse direction can be neglected. Thus, the governing differential equations for the deflection of the plate according to the Winkler model and the linear elasticity framework of the Kirchhoff-Love theory are based on the following fourth-order Dirichlet boundary value problem: 
\begin{equation}
	\begin{aligned}
		\div\left(\div\left(\frac{h^3}{12}
		\mathbb{C}\nabla^2w\right)\right)+kw & = f\delta(P_0)  &\text{in}\; & \Omega,\\
		                                   w & = 0             &\text{on}\; & \partial\Omega,\\ 
		       \frac{\partial w}{\partial n} & = 0             &\text{on}\; & \partial\Omega,
	\end{aligned}
	\label{eq:governing_equations}
\end{equation}
where $\mathbb{C}=\frac{E}{1-\nu^2}$ is the elasticity tensor, with $E$ and $\nu$ that represent the Young's modulus and the Poisson's ratio of the plate material, respectively, while $n$ is the unit outer normal to the boundary $\partial\Omega$.  

\begin{remark}
	Let us recall that the divergence of a second order tensor $\mathbb{T}=\mathbb{C}\nabla^2$, with $\mathbb{C}(\cdot)\in L^\infty(\Omega,\mathcal{L}(\mathbb{M}^2,\mathbb{M}^2))$ is defined as follows:
	$$
	(\div\mathbb{T})_i=\sum_j\partial_{x_j}(T_{ij}),
	$$
	hence, we are sure that we can reapply the divergence operator to $\div\mathbb{T}$ in the equation (\ref{eq:governing_equations}).
\end{remark}

We assume the flexural rigidity $\mathbb{D}=\frac{h^3\mathbb{C}}{12}$ to be constant, so that the principal part of the equation (\ref{eq:governing_equations}) is the biharmonic operator, then we adopt a finite difference method (FDM) for the bi-dimensional case to discretize the outgoing biharmonic operator and solve the fourth-order Dirichlet boundary value problem (\ref{eq:governing_equations}) obtaining a numerical solution of the deflection variable $w$. 
We are mainly interested in the development of an ensemble Kalman filter method for our inverse problem in order to recover successfully the Winkler coefficient.

The Ensemble Kalman inversion (EKI) is a recent technique \cite{schillings2018convergence,Chada_2018,chada2019incorporation,kovachki2019ensemble,ding2021ensemble,chada2020tikhonov} whose purpose is to solve Bayesian inverse problems through data assimilation methodologies. It is an iterative method that deals with static problems based on the ensemble Kalman filter (EnKF) which is an optimal algorithm originally designed for state estimation of dynamical systems \cite{evensen2003ensemble} and it performs well when applied to inverse problems \cite{herty2019kinetic,schillings2017analysis,iglesias2013ensemble,herty2020continuous,schwenzer2020identifying}. Although in this algorithm the means and covariances are computed from empirical ensemble, it works as an optimizer which requires the evaluation of a forward map $\mathcal{G}:\mathcal{X}\to\mathcal{Y}$ taking the unknown Winkler coefficient  $k\in\mathcal{X}$, $k>0$ from noisy observational data $y\in\mathcal{Y}$ such that:
\begin{equation}
	y=\mathcal{G}(k)+\eta,
	\label{eq:noisy_obs_data}
\end{equation}
where $\eta$ denotes the noise in the measurements. We assume that $\mathcal{X}$ and $\mathcal{Y}$ are Hilbert finite dimensional spaces, i.e. $\mathcal{X}\equiv\mathbb{R}^p_+$ and $\mathcal{Y}\equiv\mathbb{R}^q$, with $p,q\in\mathbb{N}$. The noise $\eta$ is a realization from the Gaussian random variable $\mathcal{N}(0,\Gamma)$, where the covariance $\Gamma=\gamma\, \textrm{Id}$, with $\gamma>0$ representing the noise level. The general inverse problem (\ref{eq:noisy_obs_data}) is equivalent to reconstruct $k$ from $y$, therefore we minimize the least squares objective function defined by 
\begin{equation}
	\Phi(k\,;y)=\frac{1}{2}\left\Vert y-\mathcal{G}(k)\right\Vert^2_\Gamma
\end{equation}
where we define $\Vert\,\cdot\,\Vert_\Gamma=\Vert\Gamma^{-1/2}\,\cdot\Vert_{\mathcal{Y}}$. Here, the operator $\Gamma$ acts on the measurement accuracy, so the absolute model-data misfit $(y-\mathcal{G}(k))$ can be weighted to improve the quality of the measurements. Note that we loose the uniqueness of the recovering variable $k$, because the inverse problems are ill-posed on $\mathcal{Y}$, hence minimizing the objective function $\Phi$ in $\mathcal{X}$ requires some form of regularization \cite{engl1996regularization}.

The study of inverse problems consists in the recovery of quantities of interest, slightly perturbed by measurements affected by noise. What has been done in the literature is a type of approach that attempts to deduce the unknown variables by minimizing, in a suitable norm, the difference between the measurements and the solution found through the mathematical model to be used. Another widely used approach is the probabilistic Bayesian approach \cite{stuart2010inverse, kaipio2006statistical}, which sees $k$, $y$ and $\eta$ as random variables and focuses on the probabilistic distribution $u|y$ constructed via Bayes's theorem. Some of these methods are more traditional such as Markov chain Monte Carlo (MCMC) methods \cite{bertozzi2018uncertainty} which characterize the uncertainty of results through statistical properties, and others more recent, such as the ensemble Kalman inversion (EKI) which solves Bayesian inverse problems through data assimilation methodologies \cite{stuart2015data,li2008numerical,evensen2009data}.\\

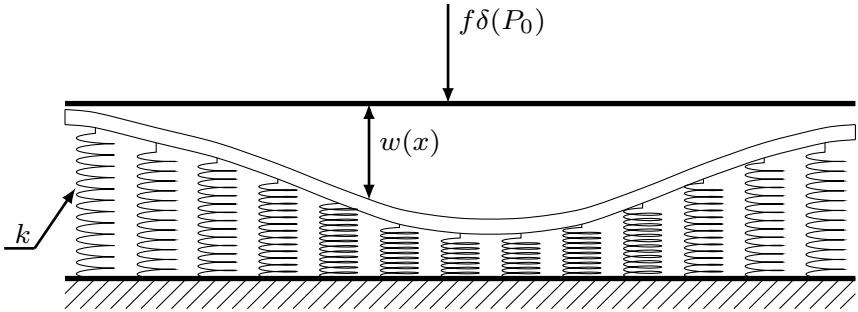
\begin{figure}[!ht]
	\centering
	\begin{tikzpicture}[scale=4]
		\draw[line width=2pt] (0,0) -- (2.6,0);
		\draw (0.05,0) -- (0,-0.05);
		\foreach \x in {0,0.05,...,2.55}
		{
			\draw (0.1+\x,0) -- (\x,-0.1);
		}
		\draw (2.6,-0.05) -- (2.55,-0.1);
		\draw [decoration={aspect=0.1, segment length=1.5mm, amplitude=2.5mm,coil},decorate] (0.1,0) -- (0.1,0.5);
		\draw [decoration={aspect=0.1, segment length=1.3mm, amplitude=2.5mm,coil},decorate] (0.3,0) -- (0.3,0.45);
		\draw [decoration={aspect=0.1, segment length=1.1mm, amplitude=2.5mm,coil},decorate] (0.5,0) -- (0.5,0.4);
		\draw [decoration={aspect=0.1, segment length=0.9mm, amplitude=2.5mm,coil},decorate] (0.7,0) -- (0.7,0.33);
		\draw [decoration={aspect=0.1, segment length=0.7mm, amplitude=2.5mm,coil},decorate] (0.9,0) -- (0.9,0.25);
		\draw [decoration={aspect=0.1, segment length=0.66mm, amplitude=2.5mm,coil},decorate] (1.1,0) -- (1.1,0.18);
		\draw [decoration={aspect=0.1, segment length=0.65mm, amplitude=2.5mm,coil},decorate] (1.3,0) -- (1.3,0.15);
		\draw [decoration={aspect=0.1, segment length=0.65mm, amplitude=2.5mm,coil},decorate] (1.5,0) -- (1.5,0.15);
		\draw [decoration={aspect=0.1, segment length=0.66mm, amplitude=2.5mm,coil},decorate] (1.7,0) -- (1.7,0.18);
		\draw [decoration={aspect=0.12, segment length=0.7mm, amplitude=2.5mm,coil},decorate] (1.9,0) -- (1.9,0.25);
		\draw [decoration={aspect=0.1, segment length=0.9mm, amplitude=2.5mm,coil},decorate] (2.1,0) -- (2.1,0.33);
		\draw [decoration={aspect=0.1, segment length=1.1mm, amplitude=2.5mm,coil},decorate] (2.3,0) -- (2.3,0.4);
		\draw [decoration={aspect=0.1, segment length=1.3mm, amplitude=2.5mm,coil},decorate] (2.5,0) -- (2.5,0.45);
		\draw  plot [smooth] coordinates {(0,0.51) (0.1,0.5) (0.3,0.45) (0.5,0.4) (0.7,0.33) (0.9,0.25) (1.1,0.18) (1.3,0.15) (1.5,0.15) (1.7,0.18) (1.9,0.25) (2.1,0.33) (2.3,0.4) (2.5,0.45) (2.6,0.46) };
		\draw (2.6,0.46) -- (2.6,0.51);
		\foreach \x in {0.05}
		{
			\draw  plot [smooth] coordinates {(0,0.51+\x) (0.1,0.5+\x) (0.3,0.45+\x) (0.5,0.4+\x) (0.7,0.33+\x) (0.9,0.25+\x) (1.1,0.18+\x) (1.3,0.15+\x) (1.5,0.15+\x) (1.7,0.18+\x) (1.9,0.25+\x) (2.1,0.33+\x) (2.3,0.4+\x) (2.5,0.45+\x) (2.6,0.46+\x)};
		}
		\draw (0,0.51) -- (0,0.56);
		
		\draw[line width=2pt] (0,0.58) -- (2.6,0.58);
		\draw[-latex,line width=1pt] (1.26,0.91) -- (1.26,0.58);
		\draw[-latex,line width=1pt] (-0.1,0.1) -- (0.035,0.3);
		\draw[line width=1pt] (-0.2,0.1) -- (-0.1,0.1);
		\draw[latex-latex, line width=1pt] (1,0.58) -- (1,0.26);
		\node [right,scale=1.3] at (1,0.45) {$w(x)$};
		\node [right,scale=1.3] at (-0.2,0.15) {$k$};
		\node [right,scale=1.3] at (1.26,0.85) {$f\delta(P_0)$};
	\end{tikzpicture}
	\caption{Equivalent foundation resting on Winkler spring bed.}
	\label{fig:Winkler_bed}
\end{figure}

We deal with a deformable plate resting on the soil. We have adopted the Winkler model, one of the most commonly used classical mathematical model to describe the behavior of the plate subject to displacement due to an applied load. In this case the Winkler model compares the settlement of the plate due to the application of a load at a specific point $P_0$ of the plate of intensity $f>0$. In Figure \ref{fig:Winkler_bed} we show the physical representation of the Winkler foundation.
For the success of the EKI method, the prior measure $\mu_0$ plays a fundamental role, since it affects the regularization of the method. For the our experiments, the initial ensemble of particles is chosen as a normal distribution $\mu_0=\mathcal{N}(0,C_0)$, with covariance operator $C_0=\beta(\Delta^2-kI)^{-1}$ related to a Brownian bridge as in \cite{schillings2017analysis}, where the factor $\beta\in\mathbb{R}$ depends on the simulation.
\subsection{Outline} The article is structured as follows: in Section \blue{\ref{section_fd}}, we introduce the finite difference method for the plate in the bi-dimensional case to discretize the biharmonic operator of the governing equations. We present, in Section \blue{\ref{section_inv_prob}}, the inverse problem and define the ensemble Kalman inversion method to obtain an optimal estimation of the subgrade coefficient. Numerical experiments are presented in Section \blue{\ref{section_num_res}}, distinguishing between the case with noise and the case without noise (noise-free). Conclusions, perspective and acknowledgments in Section \blue{\ref{section_sum_perp}}.

\section{Finite difference method for the thin plate}\label{section_fd}
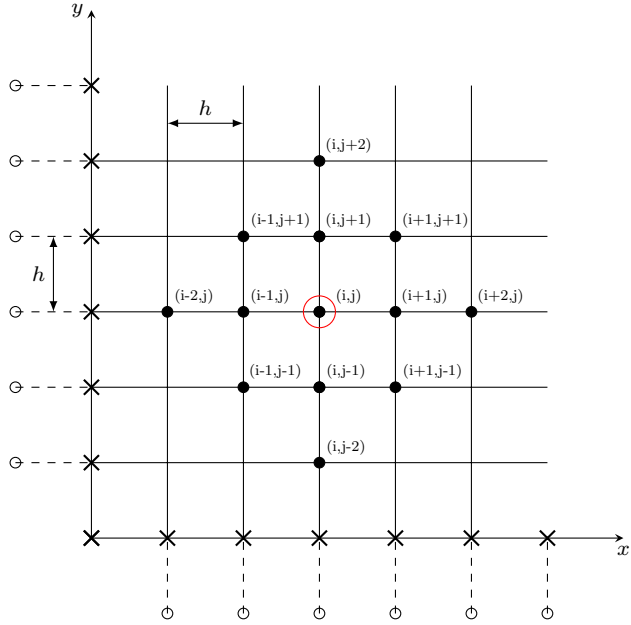
\begin{figure}[!ht]
\centering
\begin{tikzpicture}
	\draw (0, 0) grid (5, 5);
	\foreach \x in {0,1,...,5}
	{
		\draw[dashed] (-1,1+\x) -- (0,1+\x);
		\draw[dashed] (1+\x,-1) -- (1+\x,0);
		\draw (0+\x,5) -- (0+\x,6);
		\draw (5,0+\x) -- (6, 0+\x);
	}
	
	\foreach \x in {0,1,...,6}
	{
		\draw[pline] plot coordinates {(0,0+\x)};
		\draw[pline] plot coordinates {(0+\x,0)};
	}
	
	\foreach \x in {0,1,...,5}
	{
		\draw (-1,1+\x) circle (2pt);
		\draw (1+\x,-1) circle (2pt);
	}
	\draw[red] (3,3) circle (6pt);
	\foreach \x in {0,1,...,4}
	{
		\filldraw[black] (1+\x,3) circle (2pt);
		\filldraw[black] (3,1+\x) circle (2pt);
	}
	
	\foreach \x in {0,1,2}
	{
		\filldraw[black] (2+\x,4) circle (2pt);
		\filldraw[black] (2+\x,2) circle (2pt);
	}
	\draw[-stealth] (6,0) -- (7,0);
	\node[left] at (0,7) {$y$};
	\draw[-stealth] (0,6) -- (0,7);
	\node[below] at (7,0) {$x$};
	
	\node[scale=0.7, xshift=0.5cm, yshift=0.3cm] at (1,3) {(i-2,j)};
	\node[scale=0.7, xshift=0.5cm, yshift=0.3cm] at (2,3) {(i-1,j)};
	\node[scale=0.7, xshift=0.55cm, yshift=0.3cm] at (3,3) {(i,j)};
	\node[scale=0.7, xshift=0.55cm, yshift=0.3cm] at (4,3) {(i+1,j)};
	\node[scale=0.7, xshift=0.55cm, yshift=0.3cm] at (5,3) {(i+2,j)};
	\node[scale=0.7, xshift=0.5cm, yshift=0.3cm] at (3,1) {(i,j-2)};
	\node[scale=0.7, xshift=0.5cm, yshift=0.3cm] at (3,2) {(i,j-1)};
	\node[scale=0.7, xshift=0.55cm, yshift=0.3cm] at (3,4) {(i,j+1)};
	\node[scale=0.7, xshift=0.7cm, yshift=0.3cm] at (2,4) {(i-1,j+1)};
	\node[scale=0.7, xshift=0.6cm, yshift=0.3cm] at (2,2) {(i-1,j-1)};
	\node[scale=0.7, xshift=0.7cm, yshift=0.3cm] at (4,2) {(i+1,j-1)};
	\node[scale=0.7, xshift=0.75cm, yshift=0.3cm] at (4,4) {(i+1,j+1)};
	\node[scale=0.7, xshift=0.55cm, yshift=0.3cm] at (3,5) {(i,j+2)};
	
	\draw[latex-latex] (-0.5,3) -- (-0.5,4);
	\node[left] at (-0.5,3.5) {$h$};
	\draw[latex-latex] (1,5.5) -- (2,5.5);
	\node[above] at (1.5,5.5) {$h$};
	
\end{tikzpicture}
\caption{Evaluating the finite difference method at $(i,j)$.}
\label{FD_discretization}
\end{figure}

We defined the flexural rigidity of the plate as $\mathbb{D}=\frac{h^3}{12}\mathbb{C}$. For our purposes we can assume $\mathbb{D}$ to be constant, then the fourth-order Dirichlet boundary value problem (\ref{eq:governing_equations}) becomes:
\begin{equation}
\begin{aligned}
           \mathbb{D}\,\nabla^4w+kw & = f\delta(P_0) &\text{in}\; & \Omega,\\
	                            w & = 0            &\text{on}\; & \partial\Omega,\\ 
	\frac{\partial w}{\partial n} & = 0            &\text{on}\; & \partial\Omega.
\end{aligned}
\label{eq:biharmonic_sys}
\end{equation}
We apply the finite difference method for the two-dimensional model \cite{ford2014numerical}, for which we can rewrite in explicit form the system (\ref{eq:biharmonic_sys}) as follows:
\begin{equation}
\mathbb{D}\left[\frac{\partial^4w}{\partial x^4}
              +2\frac{\partial^4w}{\partial x^2\partial y^2}
              +\frac{\partial^4w}{\partial y^4}\right]+kw=f\delta(P_0),
\label{fd_2d_discretization}
\end{equation}
with $w(x,y)=0$ and $\frac{\partial w}{\partial n}=0$ on $\partial\Omega$. The notation $\frac{\partial w}{\partial n}$ refers the normal derivative, where $n$ is the outer unit vector orthogonal to the boundary $\partial\Omega$. 
The domain $\Omega$ can be partitioned into an $(n+1)\times(n+1)$ grid points with equal space step $h_s$ in the $x$ and $y$ directions. Here, we adopt a 13-point central difference formula to approximate (\ref{fd_2d_discretization}) as follows:
\begin{equation}
	\begin{aligned}
    \mathbb{D}\,&[20\,w_{ij}-8\left(w_{i-1,j}+w_{i+1,j}+w_{i,j-1}+w_{i,j+1}\right)\\
    &+2\left(w_{i-1,j-1}+w_{i+1,j-1}+w_{i-1,j+1}+w_{i+1,j+1}\right)\\
    &+\left(w_{i-2,j}+w_{i+2,j}+w_{i,j-2}+w_{i,j+2}\right)]/h_s^4+k_{ij}w_{ij}=f\delta_{ij}(P_0),
    \end{aligned}
	\label{eq:biharmonic_FD}
\end{equation}
for $1\leqslant i,j\leqslant n-1$. In Figure (\ref{FD_discretization}) we denote by crosses the nodes at the boundary where we impose $w=0$, by the filled circles the internal nodes, while by the empty circles, outside the domain, the points where $\frac{\partial w}{\partial n}=0$.

In this section, our purpose is to solve the system (\ref{eq:biharmonic_sys}) directly, by using the FDM to discretize the biharmonic operator and find the displacement $w$ given an initial value for the subgrade coefficient $k$. Thus, from (\ref{eq:biharmonic_FD}), we consider the linear model defined as:
\[
[\mathbb{D}B+kI]w=f\delta(P_0),
\]
where $B$ is an ill-conditioned positive definite $(n-1)^2\times(n-1)^2$ block pentadiagonal coefficient matrix of the form:
$$
B = \frac{1}{h_s^4}
\left[ {\begin{array}{ccccccc}
		\ddots & \ddots & \ddots &        &        &        & \\
		& \ddots & \ddots & \ddots &        &        & \\
		\ddots &        & \ddots & \ddots & \ddots &        & \\
		& B_1    & B_2    & B_3    & B_2    & B_1    & \\
		&        & \ddots & \ddots & \ddots &        & \ddots \\
		&        &        &        &        &        & \ddots \\
\end{array} } \right]
$$
with each block $B_r$, $r=1,2,3$ of size $(n-1)\times (n-1)$, in which: $B_1$ is the identity matrix, $B_2$ has the pattern $[{\begin{array}{ccc} 2,& -8, & 2\\\end{array}}]$, while $B_3$ has the pattern that changes according to the boundary conditions being considered, i.e. for $n=5$, with the boundary conditions given by (\ref{eq:biharmonic_sys}), the sub-matrix $B_3^{'}$ represents the first and the last occurrence of the full matrix $B_3$ and the sub-matrix $B_3^{''}$ the other occurrences internal to $B_3$
\[
B_3^{'}= \left[ {\begin{array}{ccccc}
		22 & -8 & 1 &0&0\\
		-8 & 21 & -8&1&0\\
		1 &-8 & 21 & -8&1 \\
		0&1 & -8 &21&-8\\
		0&0&1 & -8 &22\\
\end{array} } \right],\quad
B_3^{''}= \left[ {\begin{array}{ccccc}
		21 & -8 & 1 &0&0\\
		-8 & 20 & -8&1&0\\
		1 &-8 & 20 & -8&1 \\
		0&1 & -8 &20&-8\\
		0&0&1 & -8 &21\\
\end{array} } \right].
\]
For as regards the source term $f\delta(P_0)$, $f>0$, since we are using a FDM to discretize the biharmonic operator, we necessarily need to approximate the Dirac delta to a Gaussian distribution of mean $P_0$, within our domain $\Omega$ and a certain appropriately chosen variance (see in Section \blue{\ref{section_num_res}}).
\section{Inverse problem}\label{section_inv_prob}
By solving the system (\ref{eq:biharmonic_sys}) directly, we were able to obtain a numerical solution of the displacement $w$. Now let us consider the inverse problem, which is what interests us, where given the displacement $w$, we want to reconstruct the Winkler coefficient $k$. In order to do this, the EKI method comes to our help, which will solve the inverse problem, giving us a reconstruction of the Winkler coefficient $k$. We are interested to show how the EKI algorithm works in this case.
\subsection{Ensemble Kalman Inversion (EKI)}
EKI is a method that over the last decade has been developed as an iterative method for solving inverse problems. It works on an ensemble of particles moving them from a prior to posterior phase in pseudo-time $n$. When applied to the general inverse problem (\ref{eq:noisy_obs_data}), in which the pseudo-time-step is denoted by $dt$, then the algorithm is updated by an ensemble candidate parameter estimates $k_n=\{k_n^{(j)}\}_{j=1}^J$, which takes the following form:
\begin{align}
\label{eq:algorithm_EKI}
k^{(j)}_{n+1}&=k^{(j)}_n+\mathcal{K}_n\left(y_{n+1}^{(j)}-\mathcal{G}(k^{(j)}_n)\right)\,,\\
\label{eq:perturb_observ}
y_{n+1}^{(j)}&=y+\xi_{n+1}^{(j)},\qquad\xi_{n+1}^{(j)}\sim\mathcal{N}(0,\Sigma)\quad\text{i.i.d.}
\end{align}
\vspace{0.1cm}
for each particle $j=1,\ldots,J$, with $J\in\mathbb{N}$. The operator
\begin{equation}
\mathcal{K}_n=C^{kw}(k_n)\left(C^{ww}(k_n)+dt^{-1}\Gamma\right)^{-1}
\label{kalman_gain}
\end{equation}
is the so-called \textit{Kalman gain}, where the operators 
$C^{ww}$, $C^{kw}$ represent the empirical covariances, defined for $k=\{k^{(j)}\}_{j=1}^J$ by:
\begin{equation}
\begin{aligned}
C^\textit{ww}(k_n) &= \frac{1}{J}\sum_{j=1}^{J}\left(\mathcal{G}(k^{(j)}_n)-\overline{\mathcal{G}}_n\right)\otimes\left(\mathcal{G}(k^{(j)}_n)-\overline{\mathcal{G}}_n\right),\quad\overline{\mathcal{G}}_n=\frac{1}{J}\sum_{j=1}^{J}\mathcal{G}(k^{(j)}_n),\\
C^\textit{kw}(k_n) &= \frac{1}{J}\sum_{j=1}^{J}\left(k^{(j)}_n-\overline{k}_n\right)\otimes\left(\mathcal{G}(k^{(j)}_n)-\overline{\mathcal{G}}_n\right),\quad\overline{k}_n=\frac{1}{J}\sum_{j=1}^{J}k^{(j)}_n.
\end{aligned}
\label{empirical_covariances}
\end{equation}
where $\otimes$ denotes the tensor product for Hilbert spaces $\mathcal{H}_1$, $\mathcal{H}_2$ defined as:
$$
z_1\otimes z_2:\mathcal{H}_1\to\mathcal{H}_2,\quad\text{with}\quad z_1\otimes z_2(q)=\left<z_2,q\right>_{H_2}\cdot z_2,\quad\forall q\in\mathcal{H}_1.
$$

\noindent The measurements $y_{n+1}^{(j)}$ are subject to perturbations due to different errors. We consider the two typical choices for the covariance $\Sigma$, in which the measurements are kept unperturbed ($\Sigma=0$) and the case in which the measurements are perturbed ($\Sigma=\Gamma$), where $\xi_{n+1}^{(j)}$ represents the realizations of the noise $\eta$ in (\ref{eq:noisy_obs_data}).

In (\ref{eq:algorithm_EKI}), the $j$-th parameter of $k_{n+1}$ is the $j$-th parameter of $k_n$ plus some operator which will involve all particles applied to the difference between the data and the model evaluation to the candidate parameter $k^{(j)}_n$.
If the current parameter $k^{(j)}_n$ does not fit the model very well, which means the difference $y-\mathcal{G}(k^{(j)}_n)$ will be large, thus we are going to make a big change to the $k^{(j)}_{n+1}$ parameter. On the other hand, if we choose a parameter that fits the data exactly, we will not modify $k^{(j)}_{n+1}$ at all. 
\subsection{Algorithm}\label{sec:algorithm}
The EKI has been developed to be an optimization method and is used as a sampling method. We consider a fixed number of samples $J$ of the parameter $k=\{k^{(j)}\}^J_{j=1}$ according to the a prior measure $\mu_0=\{k^{(j)}_0\}^J_{j=1}$, defined as an initial normal distribution. This prior $\mu_0$ has a very important role for the optimization of the method, as the initial ensemble is drawn from $\mu_0$. The data-misfit $\theta$ could lead to over-fitting the solution, so then we consider the discrepancy principle as stopping criterium to overcome this problem. The iterations of the iterative ensemble method will be stopped when $\theta\leqslant\left\Vert\eta\right\Vert.$ We summarize in a pseudo-code the EKI algorithm step by step in Algorithm \ref{alg:cap}.

The goal of this paper is to implement the EKI algorithm for solving the inverse problem (\ref{eq:noisy_obs_data}) and reconstruct numerically the subgrade coefficient of a plate in according to the Winkler model. We are going to show the numerical simulations in the bi-dimensional case of the numerical displacement of the plate, before and after a small perturbation due to the Gaussian noise in the inverse problem and the reconstructed solution of the system (\ref{eq:governing_equations}). Finally we can show an estimation of the Winkler coefficient taking into account different test functions.

\begin{algorithm}
	\caption{Ensemble Kalman Inversion}\label{alg:cap}
	\textbf{Require:}	
	\begin{algorithmic}[1]
		\State\text{Input:} ensemble particles $J>1$; pseudo-time-step $dt\geqslant1$; maximum number of iterations $N$; $\Gamma=\gamma\,\textrm{Id}$, with $\gamma>0$ as noise level; and data $y$.
		\State\text{Initial:} $\{k^{(j)}_0\}$ sampled from the initial distribution $\mu_0$.
	\end{algorithmic}
	\textbf{Run:}
	\begin{algorithmic}[1]
		\While{$n<N$}
		\State\textbf{compute} the deviation of each ensemble $j=1,\ldots,J$ from the mean $\overline{k}$
		$$
		e^{(j)}=\frac{\left\Vert k^{(j)}_{\text{EKI}}-\overline{k}\right\Vert_{L^2(\Omega)}}{\left\Vert\overline{k}\right\Vert_{L^2(\Omega)}}
		$$
		\State\textbf{compute} the residual for each ensemble $j=1,\ldots,J$ from the truth $k^\dagger$
		$$
		r^{(j)}=\frac{\left\Vert k^{(j)}_{\text{EKI}}-k^\dagger\right\Vert_{L^2(\Omega)}}{\left\Vert k^\dagger\right\Vert_{L^2(\Omega)}}
		$$
		\State\textbf{compute} the misfit for each ensemble $j=1,\ldots,J$
		$$
		\theta^{(j)}=\left\Vert y-\overline{\mathcal{G}}(k^{(j)}_{\text{EKI}})\right\Vert_\Gamma
		$$
		\If{$\theta\leqslant\Vert\eta\Vert$}
		\State\textbf{break}
		\EndIf
		\State\textbf{compute} the empirical means and covariances defined in (\ref{empirical_covariances})
		\For{1 \textbf{to} J}
		\State perturb data with $\xi_{n+1}^{(j)}$ drawn i.i.d. from (\ref{eq:perturb_observ})
	    \State\textbf{update} set $n$ to $n+1$ the Kalman gain $\mathcal{K}_n$ from (\ref{kalman_gain})
		\State\textbf{update} set $n$ to $n+1$ the iterative EKI (\ref{eq:algorithm_EKI})
		\EndFor
		\State\textbf{calculate} the mean $\overline{k}$
		\EndWhile
	\end{algorithmic}
\end{algorithm}

\section{Numerical results}\label{section_num_res}
In all the simulations we considered only the internal points, without the boundary conditions. As the inverse problem to be solved is ill-posed and physically it is not possible for us to reconstruct the conditions at the boundary. So we will show only the internal reconstruction of the related domain. In our experiments, we set the following parameters:
\begin{itemize}
\item[-] the numerical domain $\Omega=[0,1]\times[0,1]$;
\item[-] the grid points $N_x=10$, $N_y=10$;
\item[-] the pseudo-time-step $dt=1$;
\item[-] the maximum number of iterations $N=2000$;
\item[-] the point where the load is applied $P_0=(0.5,0.5)$;
\item[-] the flexural rigidity $\mathbb{D}=\textrm{Id}$;
\item[-] the concentrated force $f\delta(P_0)$, with intensity $f=1$;
\item[-] the approximation of the Dirac delta function $\delta(P_0)$ to a normal distribution $\mathcal{N}(P_0,\sigma)$, where $\sigma$ is a certain covariance matrix $$\sigma=\begin{bmatrix}
	1/s & 0\\
	0 & 1/s
\end{bmatrix},$$ with $s=10^5$.
\end{itemize}
It is natural to set in a suitable way the Young's modulus, the plate thickness and the Poisson's ratio depending on the type of the material needed and consequently set the EKI method in order to be able to reconstruct in the best way the subgrade reaction coefficient $k$.
\subsection{Test case $1$}
Here, we propose some numerical results for the model (\ref{eq:governing_equations}), useful for engineering applications. We rewrite the governing equation given in (\ref{eq:biharmonic_sys}):
$\mathbb{D}\,\nabla^4w+kw = f\delta(P_0)$ defined in $\Omega$. In this test we would like to reconstruct the Winkler coefficient $k=e^{x+y}\in L^\infty(\Omega)$.\\

\begin{figure}[!ht]
	\subfloat[]{
	\includegraphics[width=.32\textwidth]{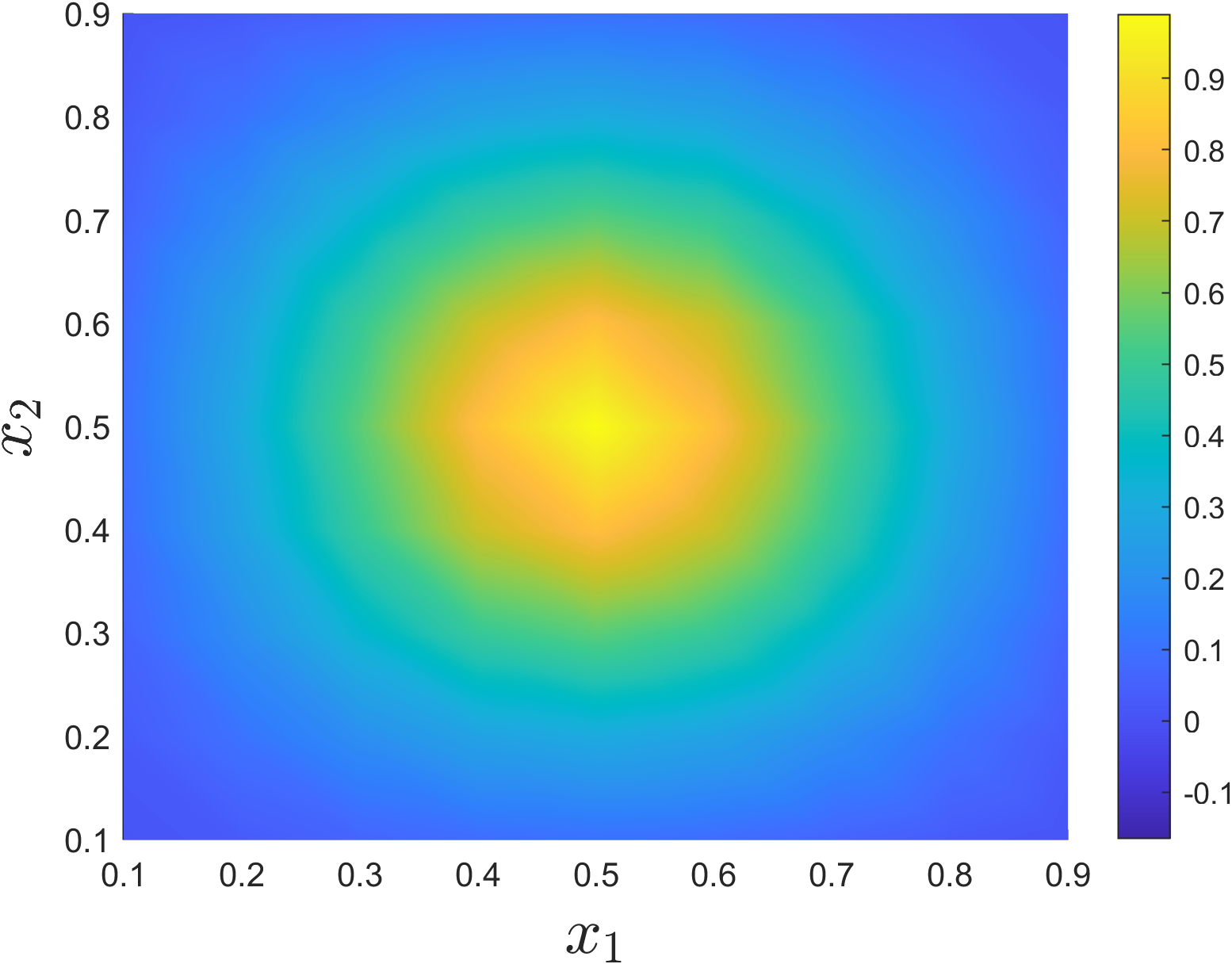}}\hspace{0.05cm}
	\subfloat[]{
	\includegraphics[width=.32\textwidth]{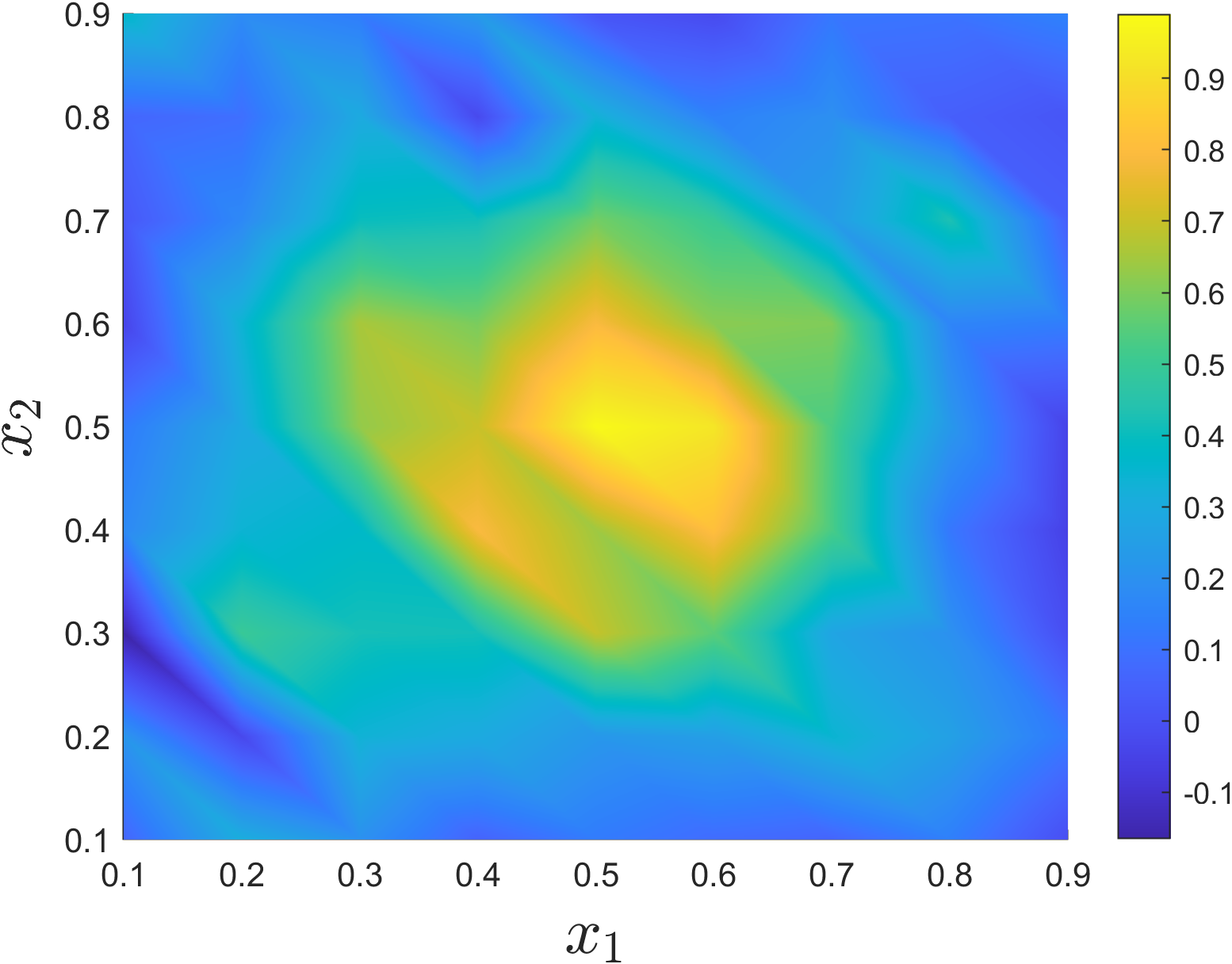}}\hspace{0.05cm}
	\subfloat[]{\includegraphics[width=.32\textwidth]{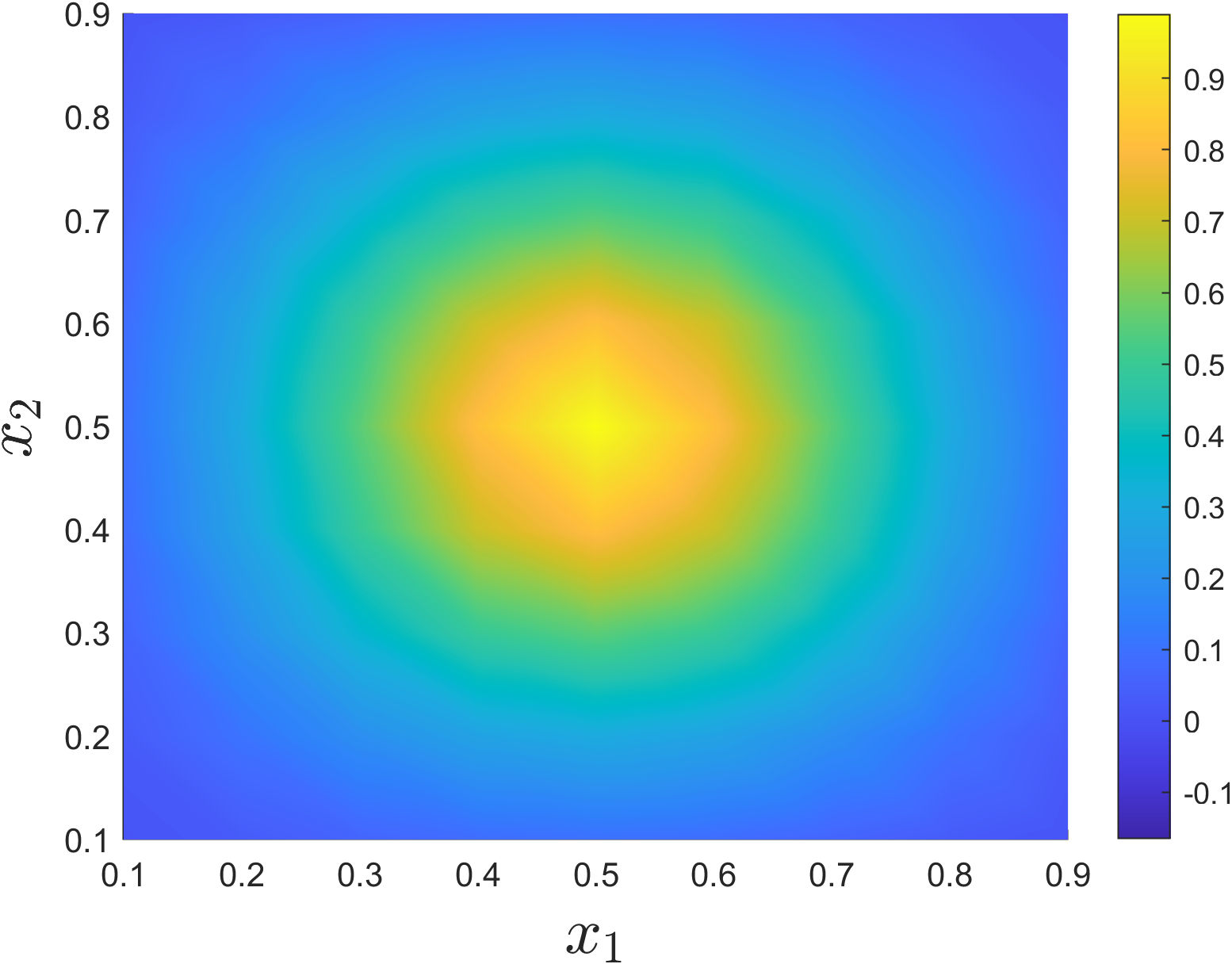}}\\
	\caption{We set the EKI algorithm by $J=1000$ particles, $\gamma=0.01$ and $\beta=10^6$. (a): Numerical solution of the displacement $w$. (b): Measurement obtained perturbing with the Gaussian noise $\eta$ the numerical solution $w$. (c): Reconstruction of the displacement $w$.}
	\label{fig:direct_problem}
\end{figure}

\noindent\textbf{The direct problem:} First of all we solve numerically the system (\ref{eq:biharmonic_sys}), directly, using the FDM to discretize the biharmonic operator, once we have the numerical solution of the system, see fig. \ref{fig:direct_problem} (a), we add a slight perturbation to the observations, adding a white noise, see fig. \ref{fig:direct_problem} (b). Finally we obtain a reconstruct of the displacement $w$ of the direct system using the EKI method, see fig. \ref{fig:direct_problem} (c).\\

\begin{figure}[!ht]
	\centering
	\subfloat[]{\includegraphics[width=.49\textwidth]{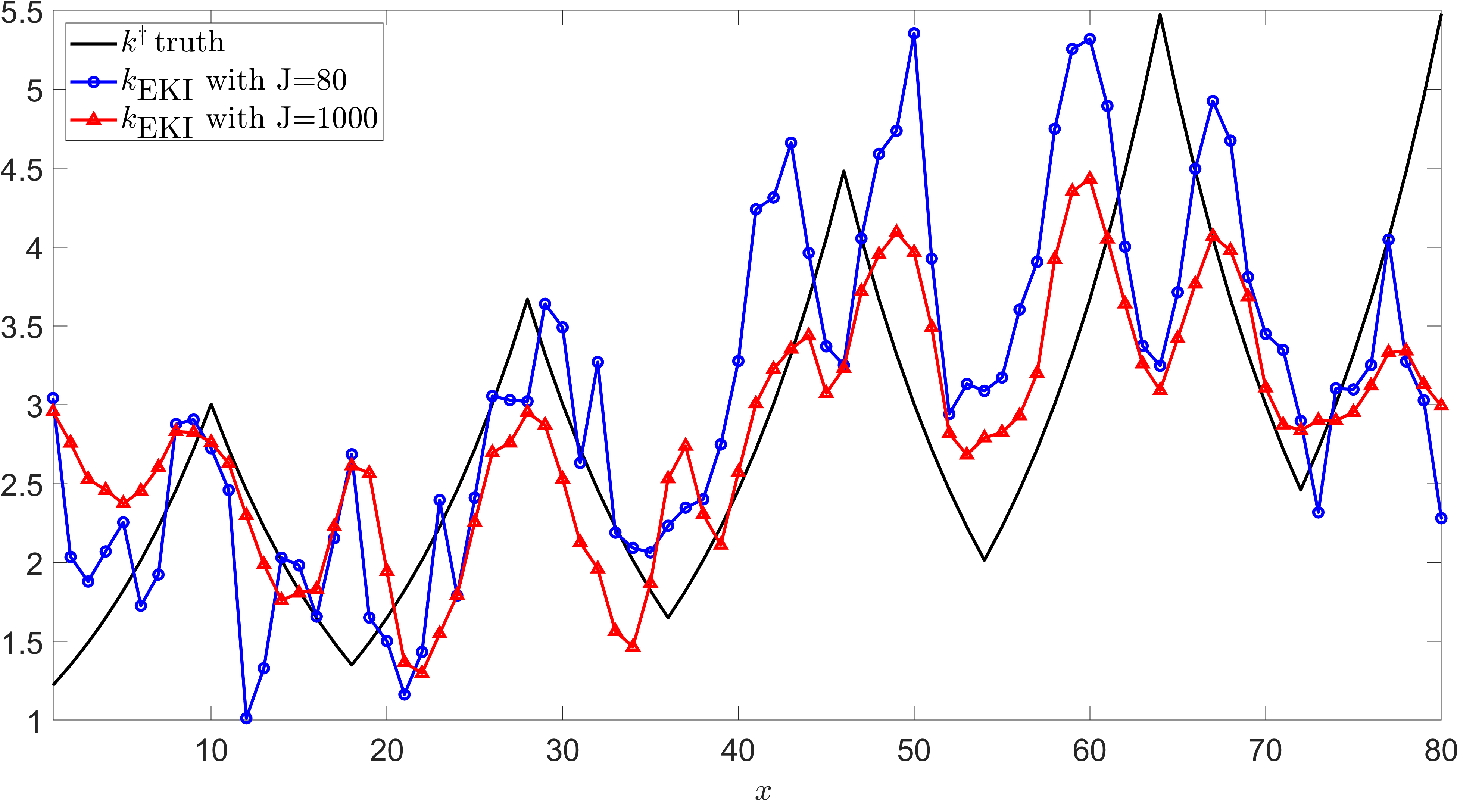}\hspace{0.1cm}}
	\subfloat[]{\includegraphics[width=.49\textwidth]{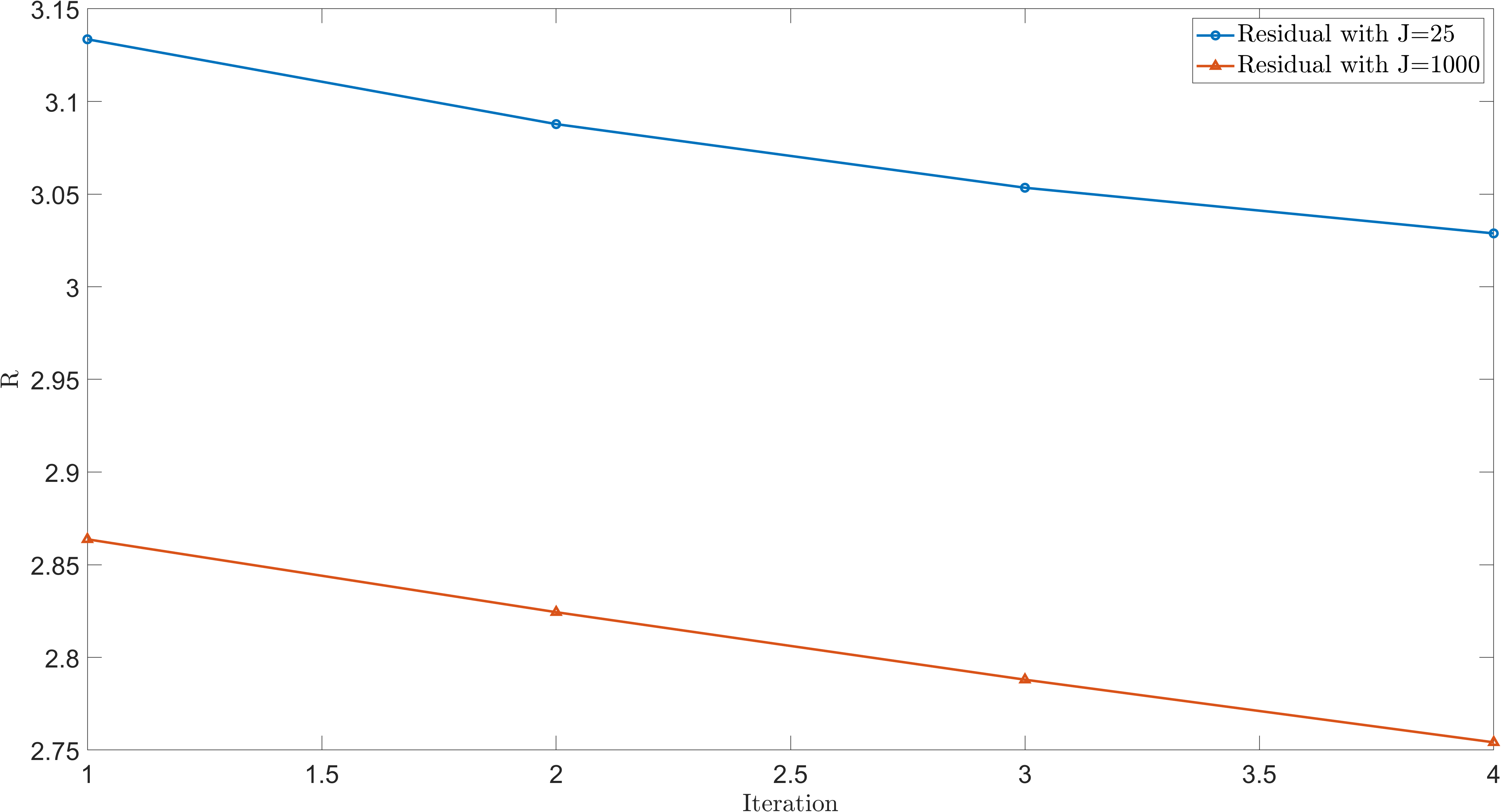}}\\
	\subfloat[]{\includegraphics[width=.49\textwidth]{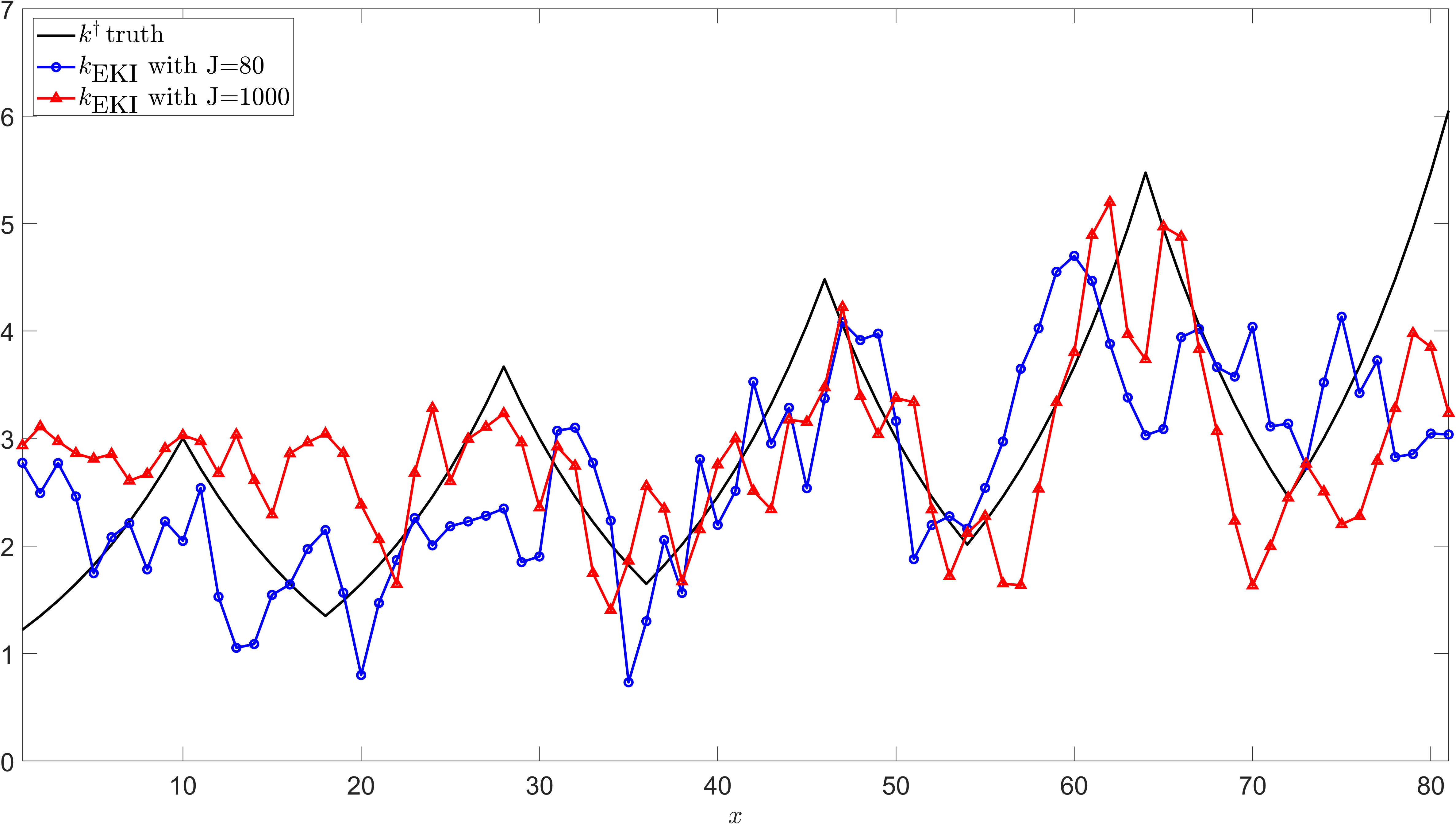}}
	\subfloat[]{\includegraphics[width=.49\textwidth]{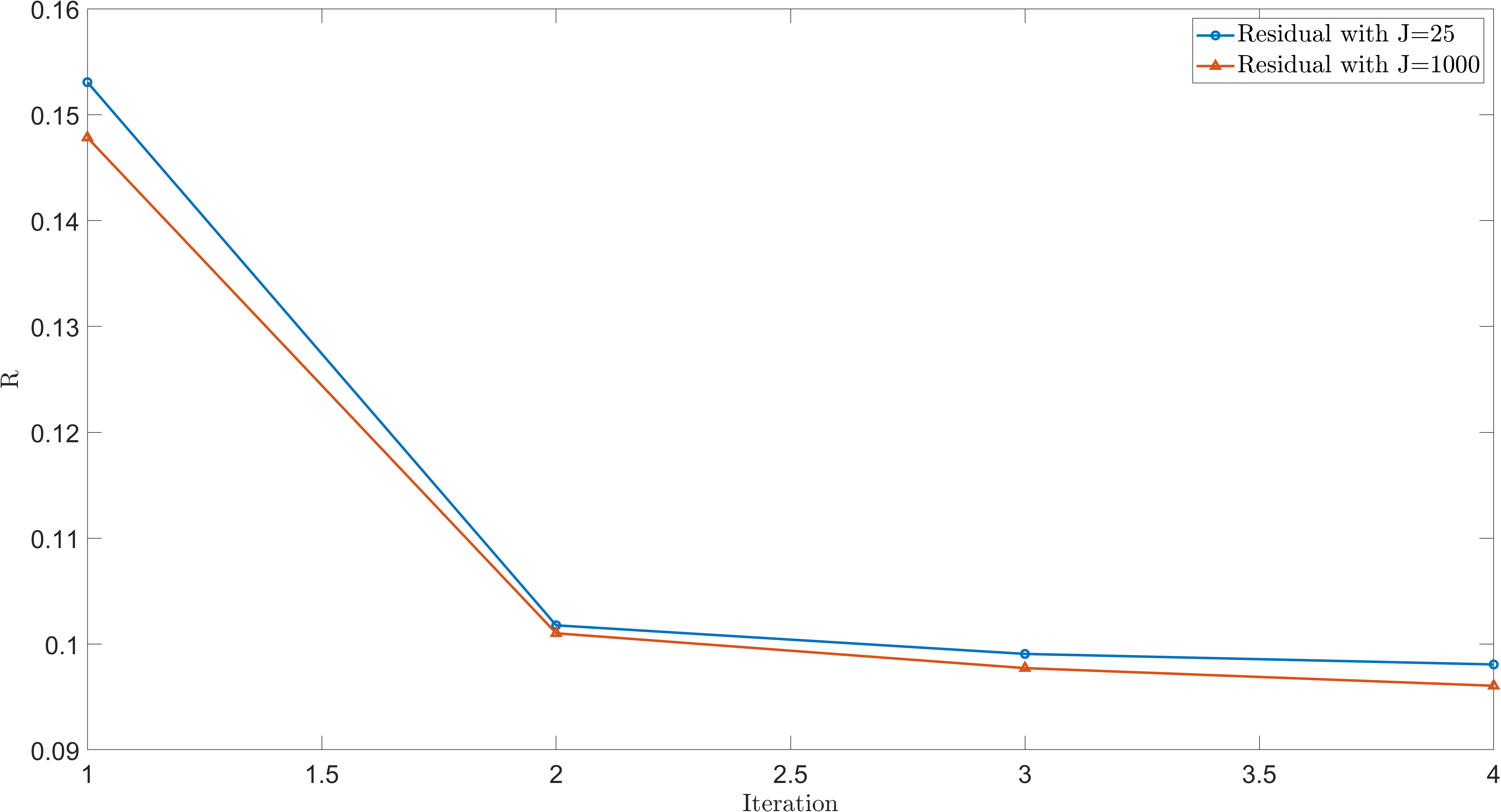}}\\
	\subfloat[]{\includegraphics[width=.49\textwidth]{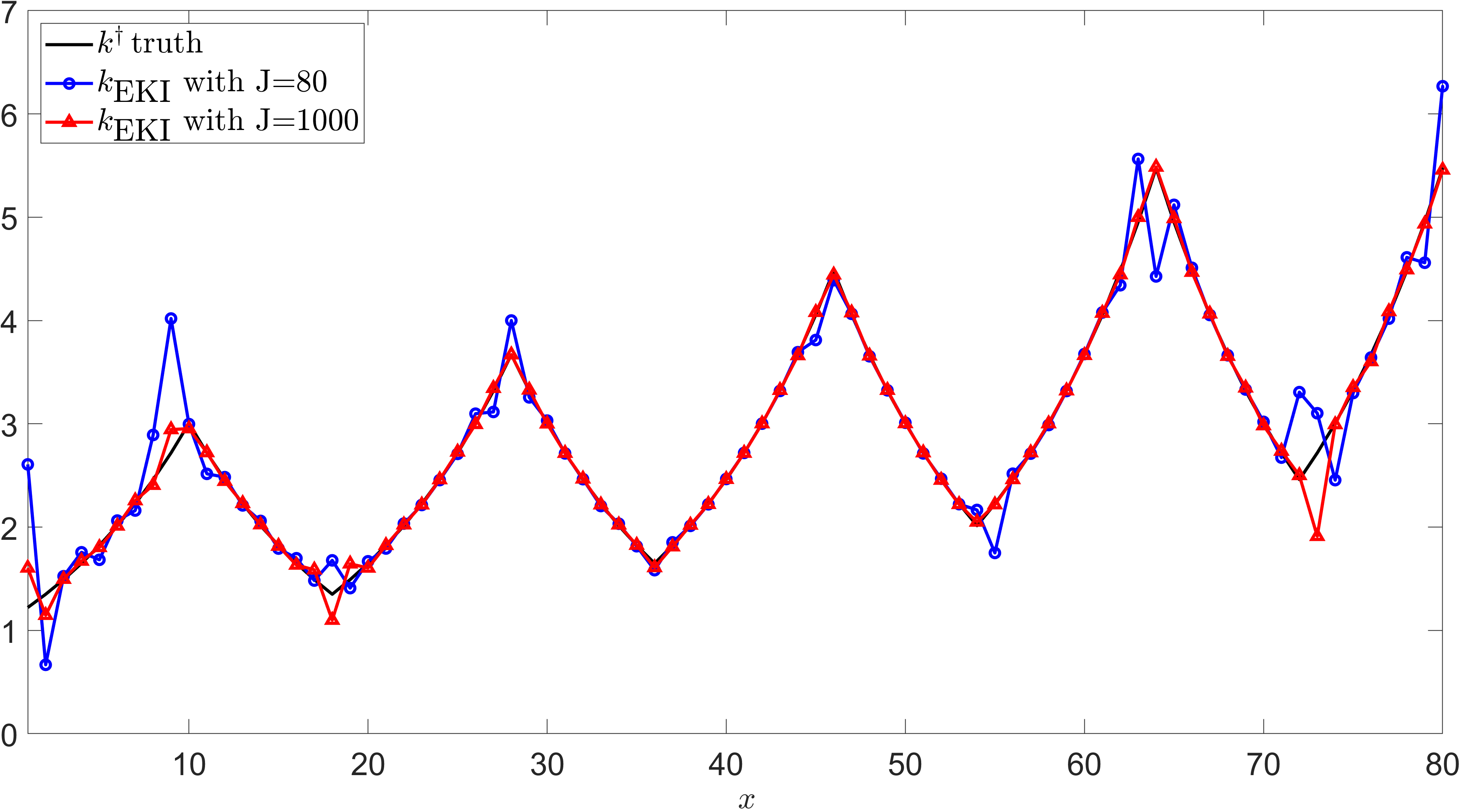}}
	\subfloat[]{\includegraphics[width=.49\textwidth]{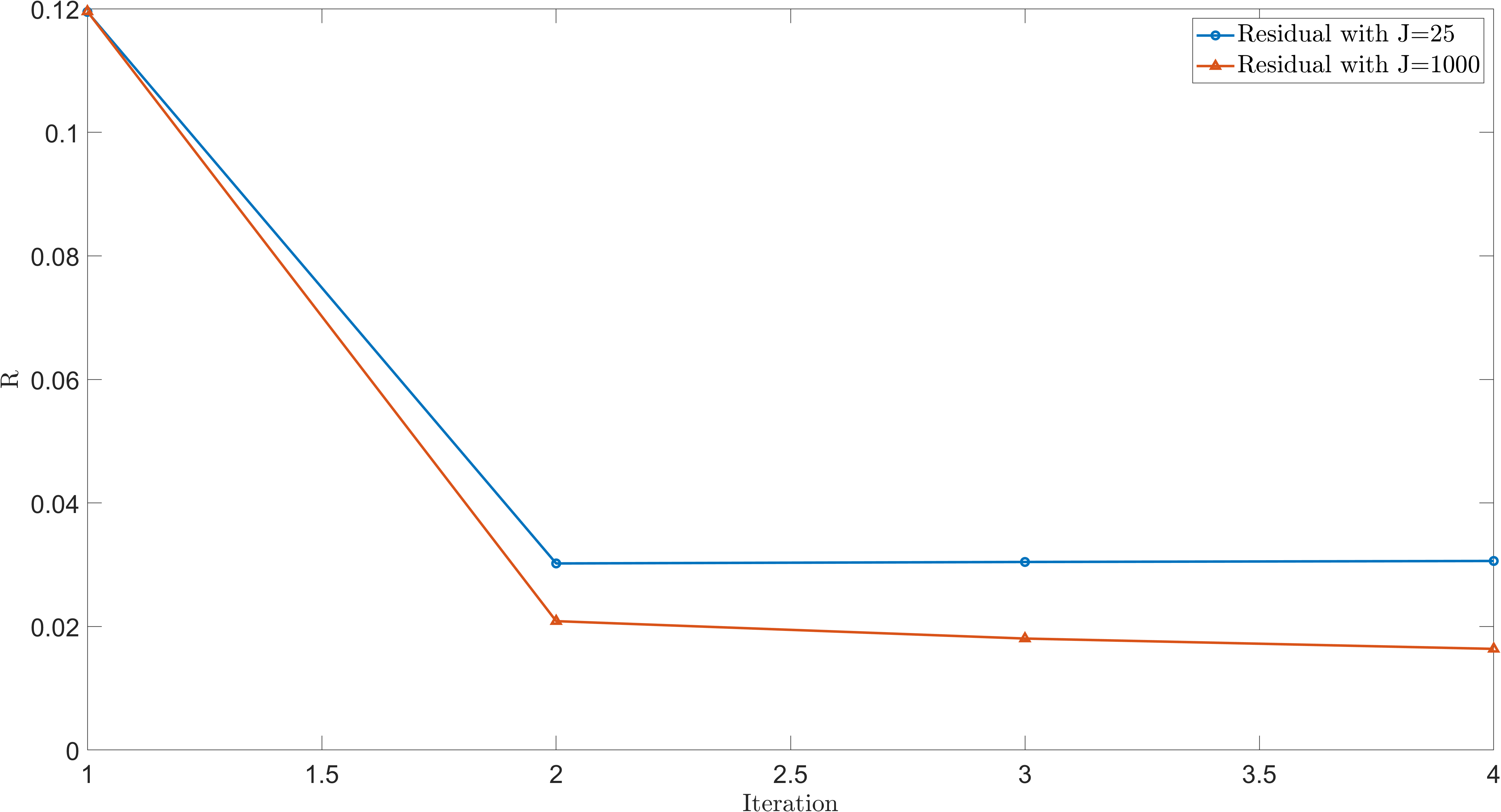}}\\
	\caption{\textit{Top line} (a),(b): with $\gamma=0.01$. \textit{Middle line} (c),(d): with $\gamma=10^{-8}$. \textit{Bottom line} (e),(f): with noise-free. \textit{Left column} (a),(c),(e): comparison between the truth value of Winkler coefficient $k^\dagger$ and its reconstruction $k_{\textrm{EKI}}$ with $J=80,1000$ particles. \textit{Right column} (b),(d),(f): residual of the Winkler coefficient $k$ for $J=25,1000$.}
	\label{residual_exp}
\end{figure}

\noindent\textbf{The inverse problem:}
Now let is analyze the inverse problem, that means, suppose we have the solution $w$ of the system (\ref{eq:biharmonic_sys}), we want to reconstruct the Winkler coefficient $k$, showing the residual that corresponds to the difference between the truth value $k^\dagger$ and the coefficient reconstructed by the EKI method $k_{\textrm{EKI}}$. We distinguish the case with noise and the case without noise and we observe that the EKI method in the absence of noise is easily able to reconstruct the Winkler coefficient $k$, while in the case in which the noise is considerable, i.e. $\gamma=0.01$, it finds some difficulty, but the residual, see the second column in fig. \ref{residual_exp}, decreases in both cases, so the EKI method is in any case a good estimator.
\begin{figure}[!ht]
	\centering
	\subfloat[]{
		\includegraphics[width=.32\textwidth]{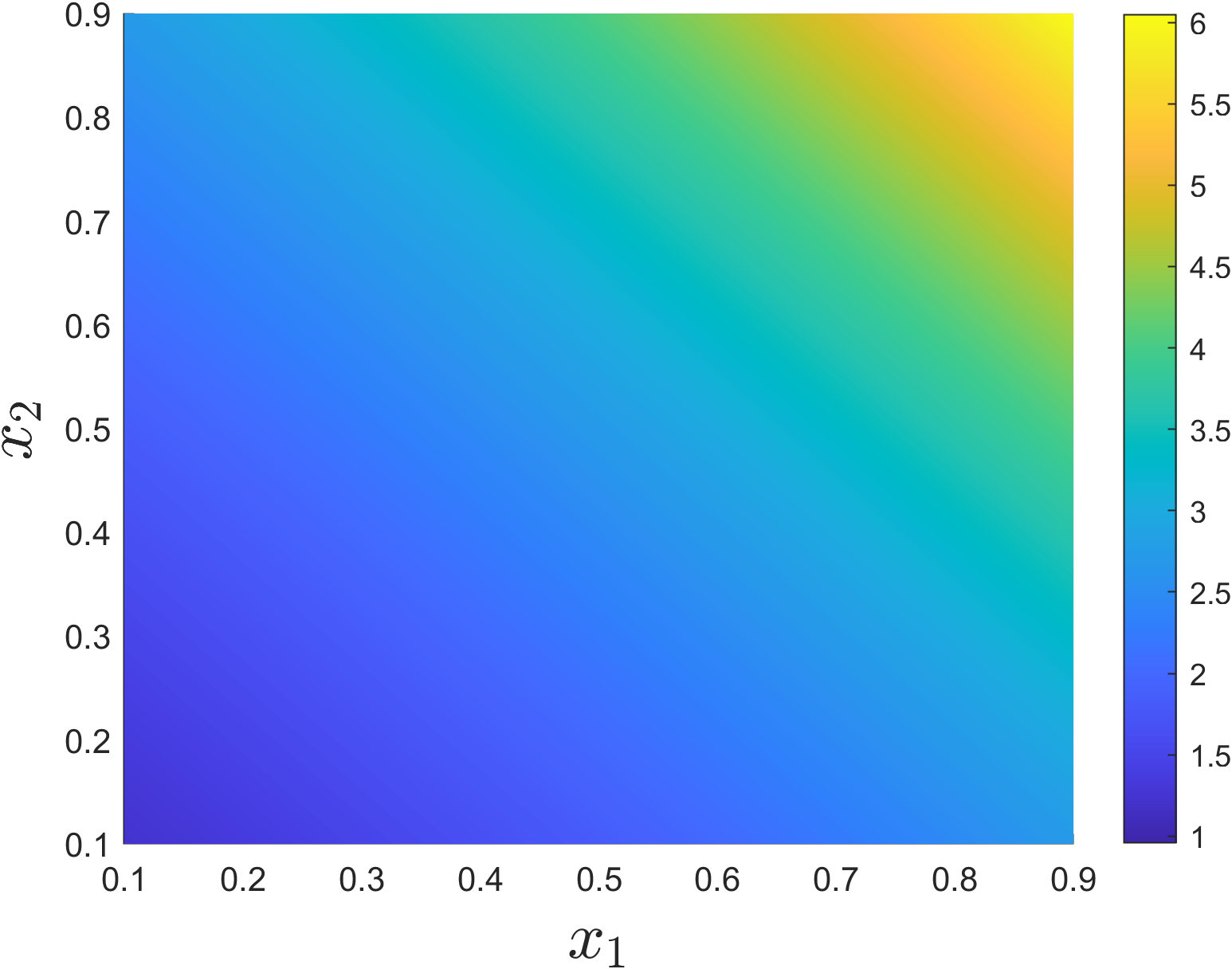}}\hspace{0.05cm}
	\subfloat[]{
		\includegraphics[width=.32\textwidth]{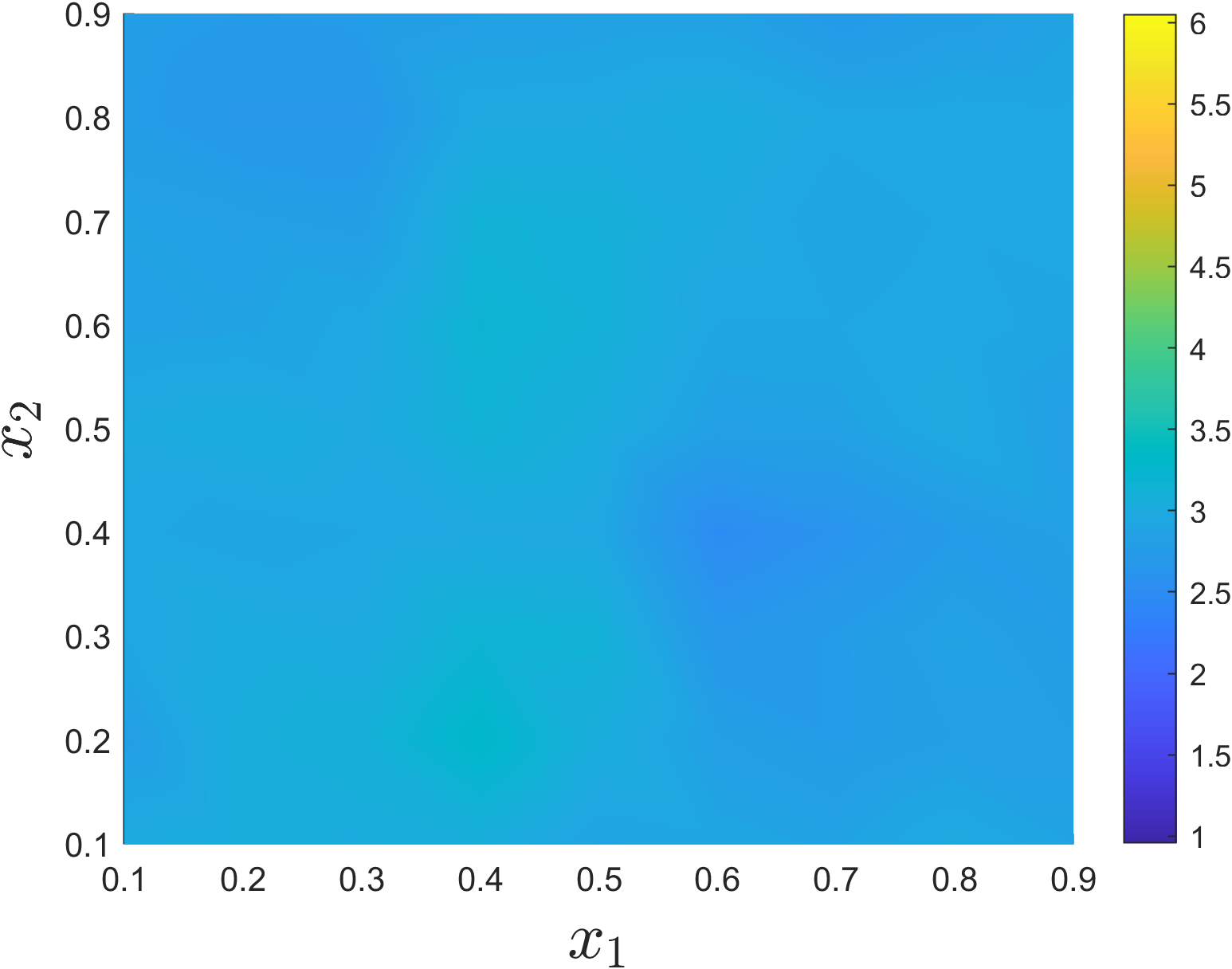}}\hspace{0.05cm}
	\subfloat[]{\includegraphics[width=.32\textwidth]{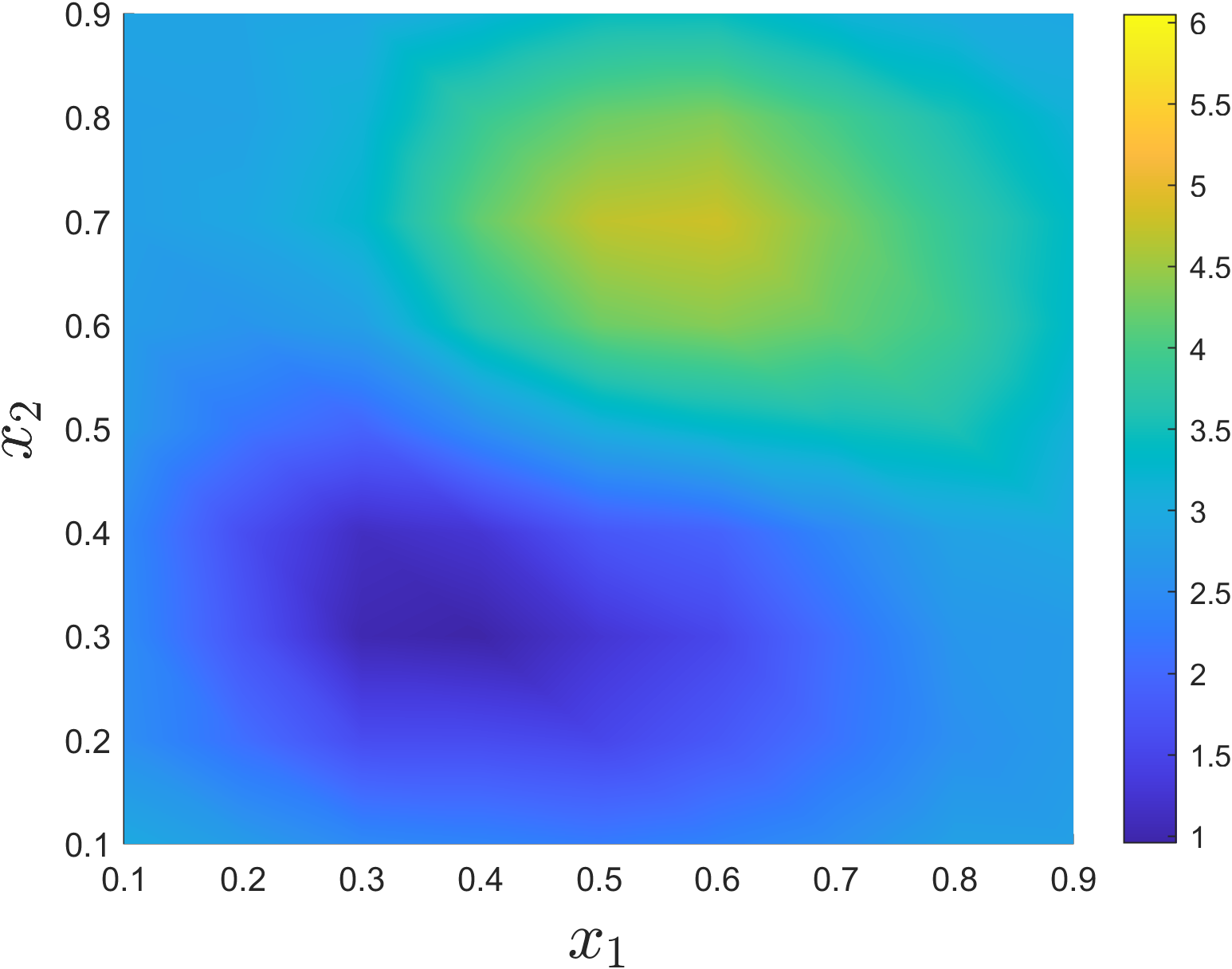}}\\
	\subfloat[]{
		\includegraphics[width=.32\textwidth]{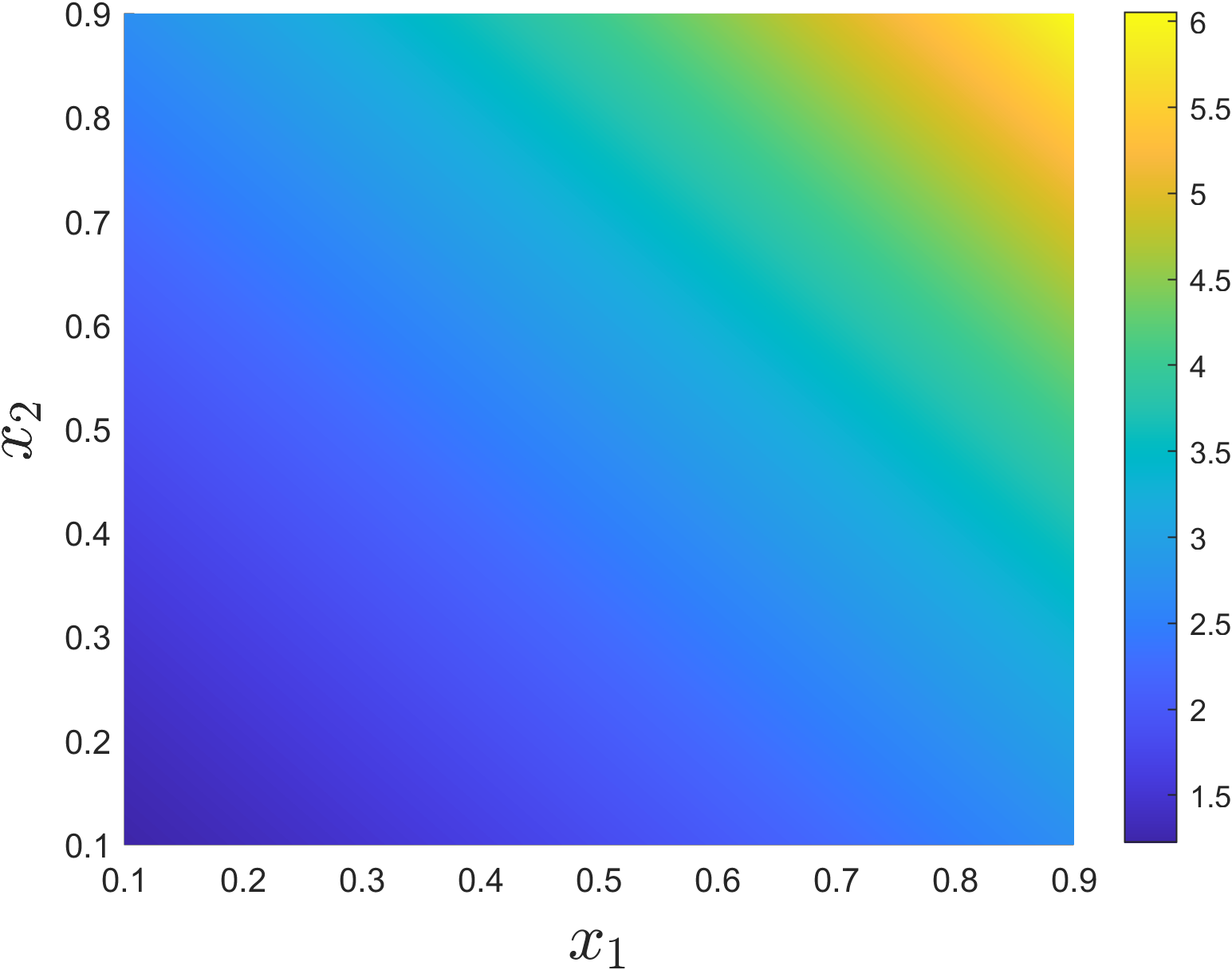}}\hspace{0.05cm}
	\subfloat[]{
		\includegraphics[width=.32\textwidth]{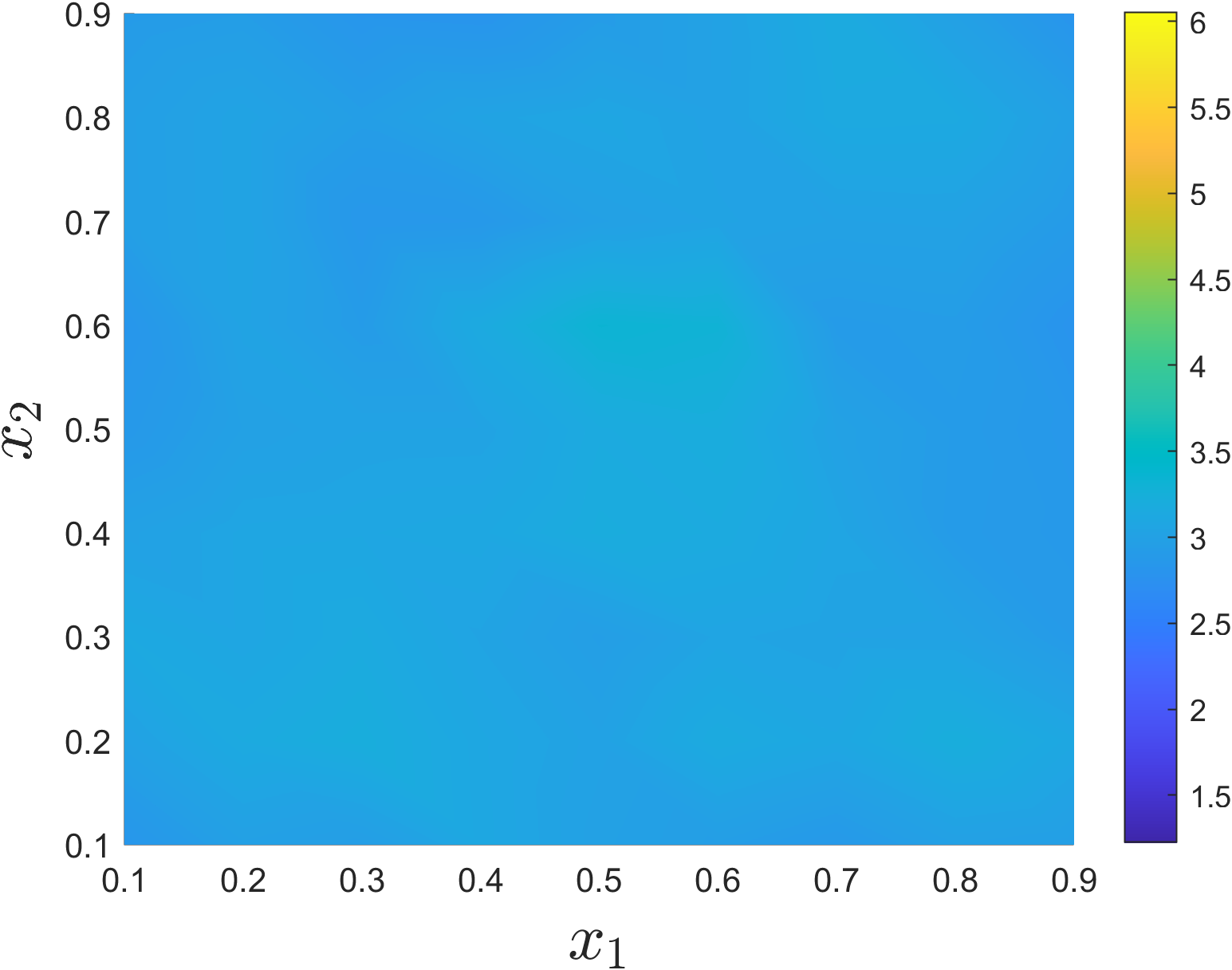}}\hspace{0.05cm}
	\subfloat[]{\includegraphics[width=.32\textwidth]{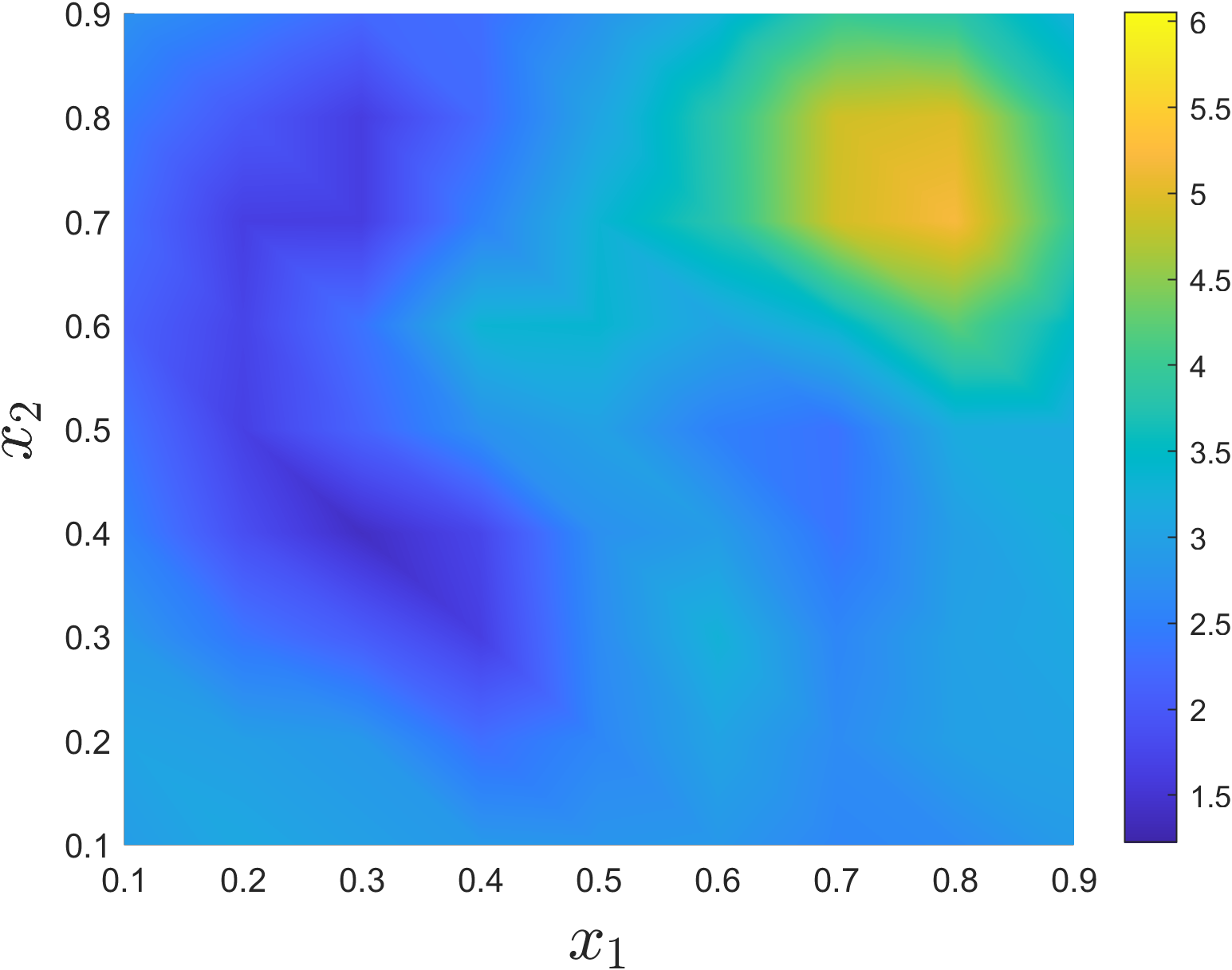}}\\
	\subfloat[]{
		\includegraphics[width=.32\textwidth]{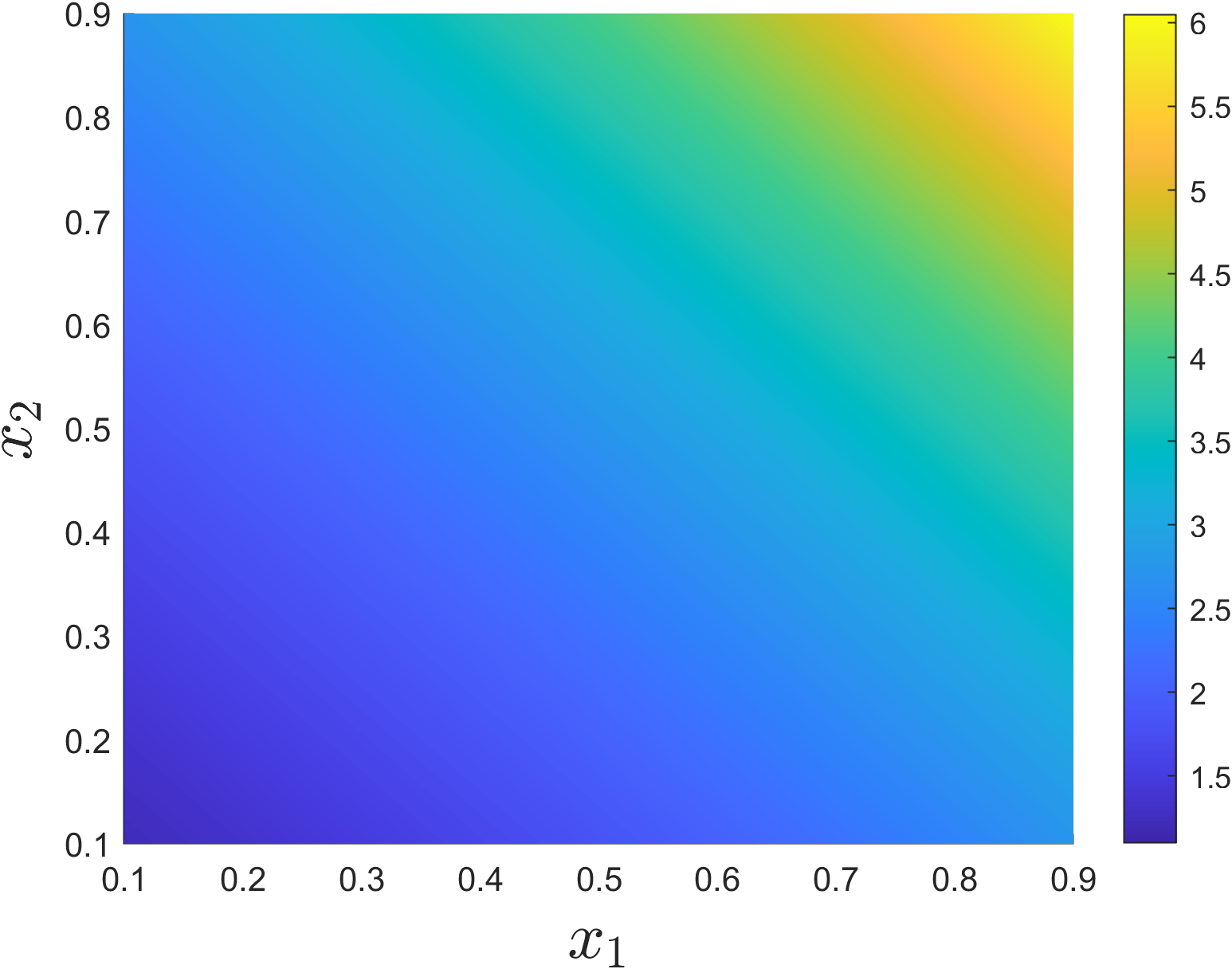}}\hspace{0.05cm}
	\subfloat[]{
		\includegraphics[width=.32\textwidth]{images/exp_J1e3_kstart_2D_noise_1em2.png}}\hspace{0.05cm}
	\subfloat[]{\includegraphics[width=.32\textwidth]{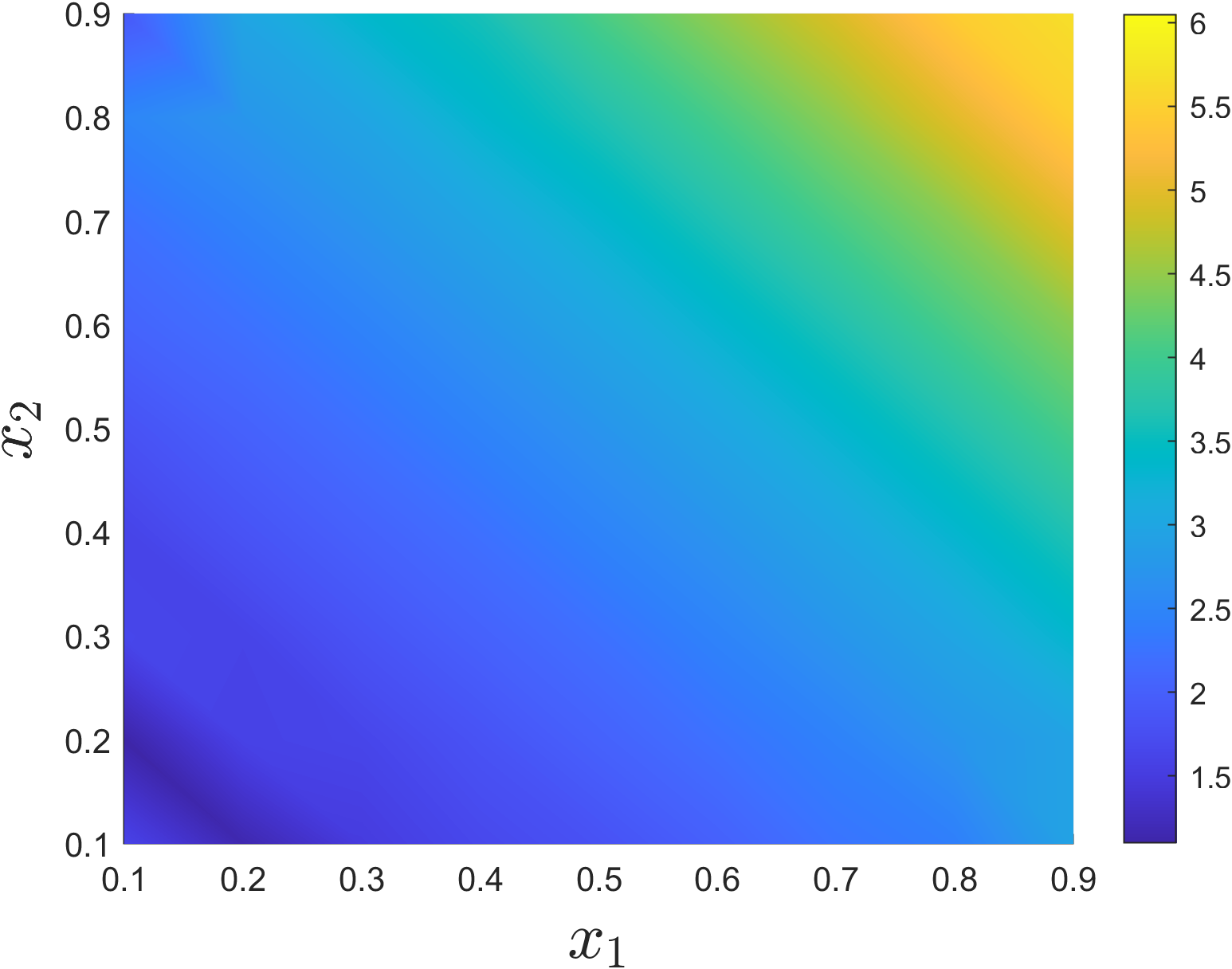}}\\
	\caption{We set $J=1000$ number of particles, $\beta=10^6$. \textit{Top line} (a),(b),(c): with $\gamma=0.01$. \textit{Middle line} (d),(e),(f): with $\gamma=10^{-8}$. \textit{Bottom line} (g),(h),(i): with noise-free. \textit{Left column} (a),(d),(g): Numerical solution of $k$. \textit{Middle column} (b),(e),(h): prior measure $\mu_0$ of $k$. \textit{Right column} (c),(f),(i): Reconstruction of $k$.}
\end{figure}
\noindent To better observe the behavior of the EKI algorithm applied to this case, we will try to plot, in the first column of the Fig. \ref{residual_exp}, in one-dimension taking the elements of the input matrix according to a snake pattern.

\subsection{Test case $2$}
Here, we consider another example for the model (\ref{eq:governing_equations}), where this time the Winkler coefficient will be defined as a piecewise constant function, since in the engineering case the subgrade reaction coefficient is a well defined number that depends on the material under examination, then, we will define a function that assumes some constant values, taken within the medium dependent range. In this test we would like to reconstruct the Winkler coefficient $k\in[0,0.15]$, defined as:
\begin{equation}
	k=
	\left\{ \begin{array}{lllll}
		0.13\quad &x\in\;]0.1,1[\;,\;& y\in\;]0,0.1[\;, \\
		0.07 &x\in\;]0.5,1[\;,& y\in[0.1,0.3[\;, \\		
		0.05 &x\in\;]0,0.9[\;,& y\in[0.3,0.5[\;, \\
		0.15 &x\in\;]0,0.6[\;,& y\in[0.5,0.7[\;, \\
		0.1 &x\in\;]0,1[\;,& y\in[0.7,1[\;.		
	\end{array} \right.
\end{equation}

\begin{figure}[!ht]
	\subfloat[]{
		\includegraphics[width=.32\textwidth]{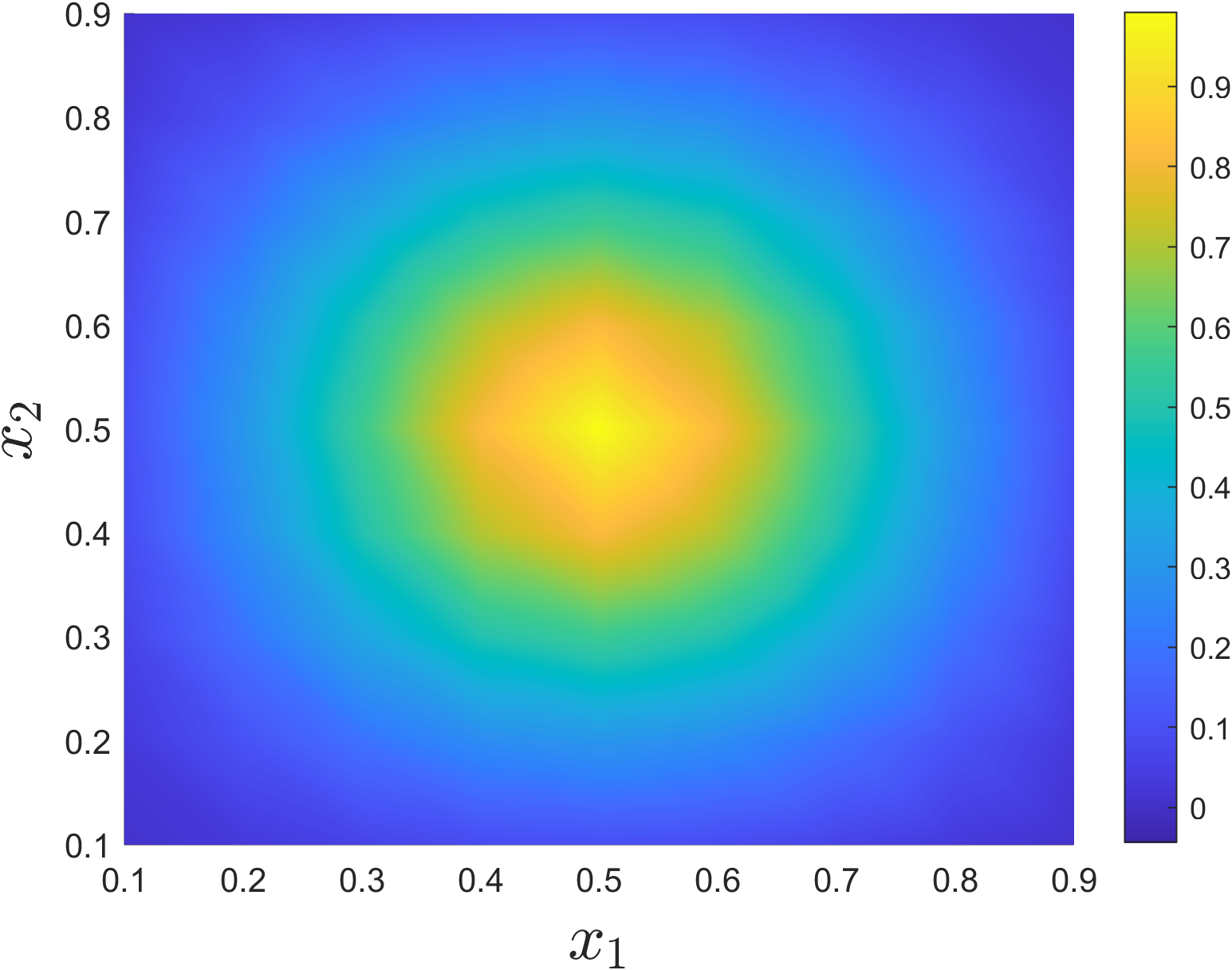}}\hspace{0.05cm}
	\subfloat[]{
		\includegraphics[width=.32\textwidth]{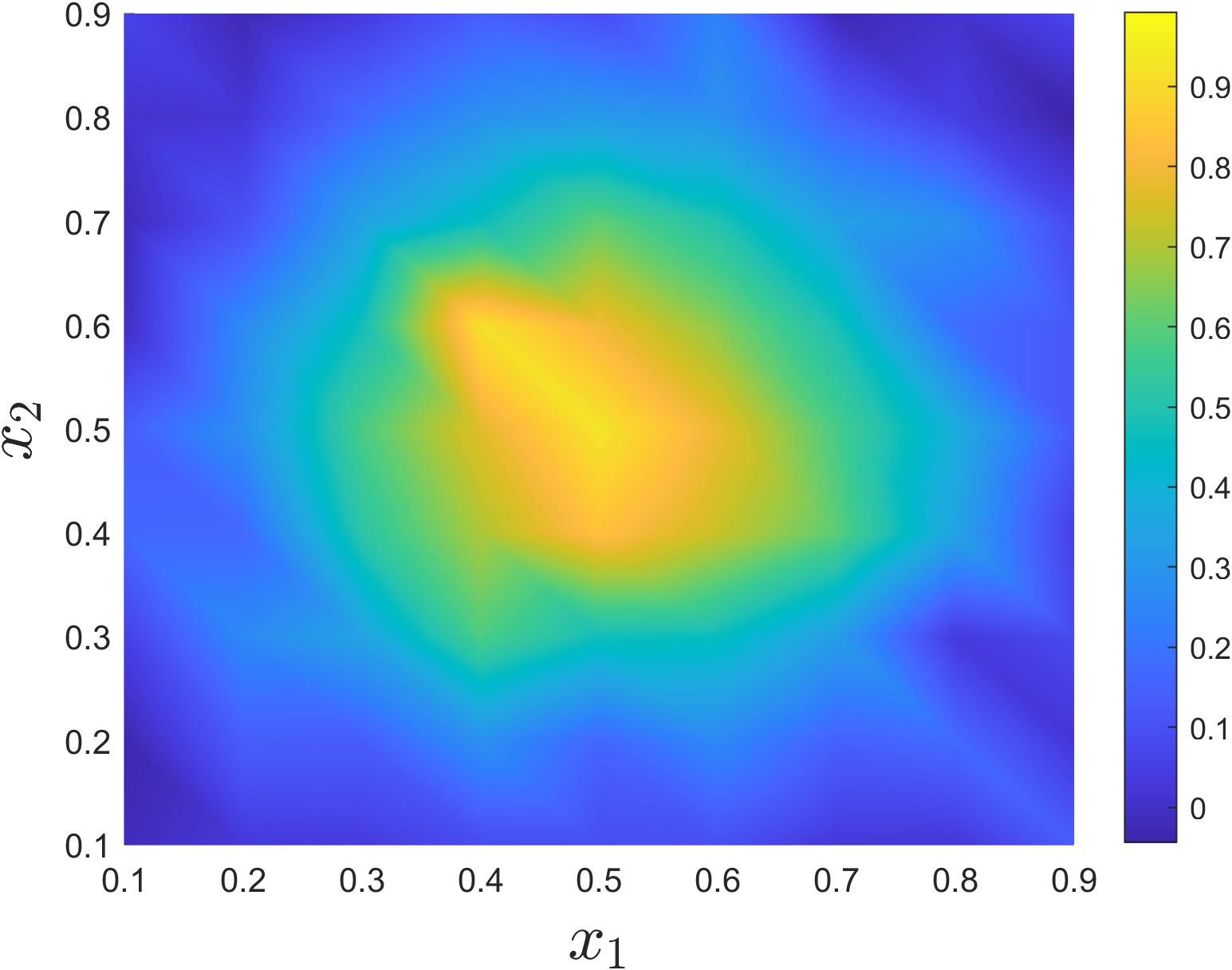}}\hspace{0.05cm}
	\subfloat[]{\includegraphics[width=.32\textwidth]{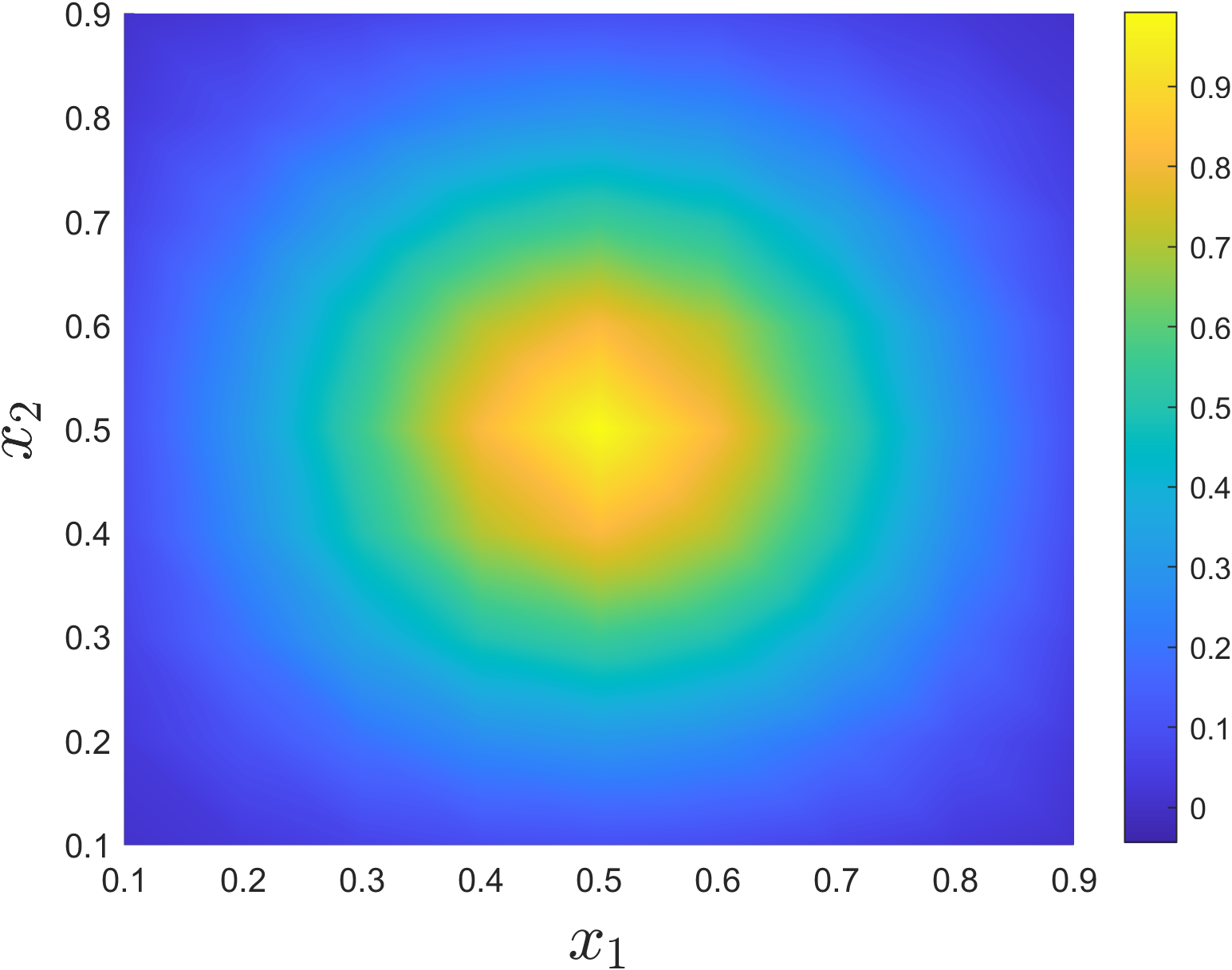}}\\
	\caption{We set the EKI algorithm by $J=1000$ particles, $\gamma=0.005$ and $\beta=6000$. (a): Numerical solution of the displacement $w$. (b): Measurement obtained perturbing with the Gaussian noise $\eta$ the numerical solution $w$. (c): Reconstruction of the displacement $w$.}
	\label{fig:direct_problem_piecewise}
\end{figure}

\noindent\textbf{The direct problem:} First of all we solve numerically the system (\ref{eq:biharmonic_sys}), directly, using the FDM to discretize the biharmonic operator, once we have the numerical solution of the system, see fig. \ref{fig:direct_problem_piecewise} (a), we add a slight perturbation to the observations, adding a white noise, see fig. \ref{fig:direct_problem_piecewise} (b). Finally we obtain a reconstruct of the displacement $w$ of the direct system using the EKI method, see fig. \ref{fig:direct_problem_piecewise} (c).\\

\begin{figure}[!ht]
	\centering
	\subfloat[]{\includegraphics[width=.49\textwidth]{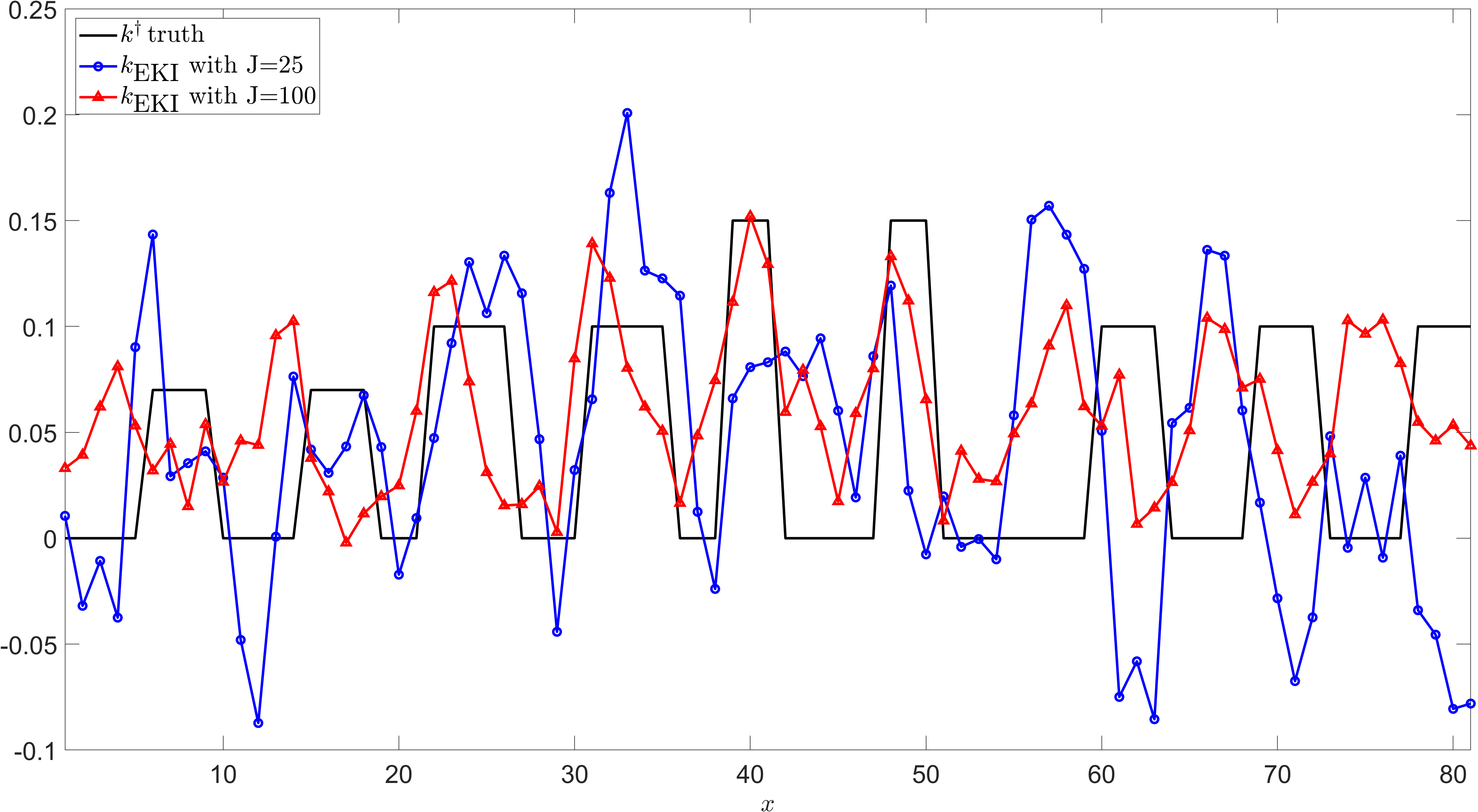}\hspace{0.1cm}}
	\subfloat[]{\includegraphics[width=.49\textwidth]{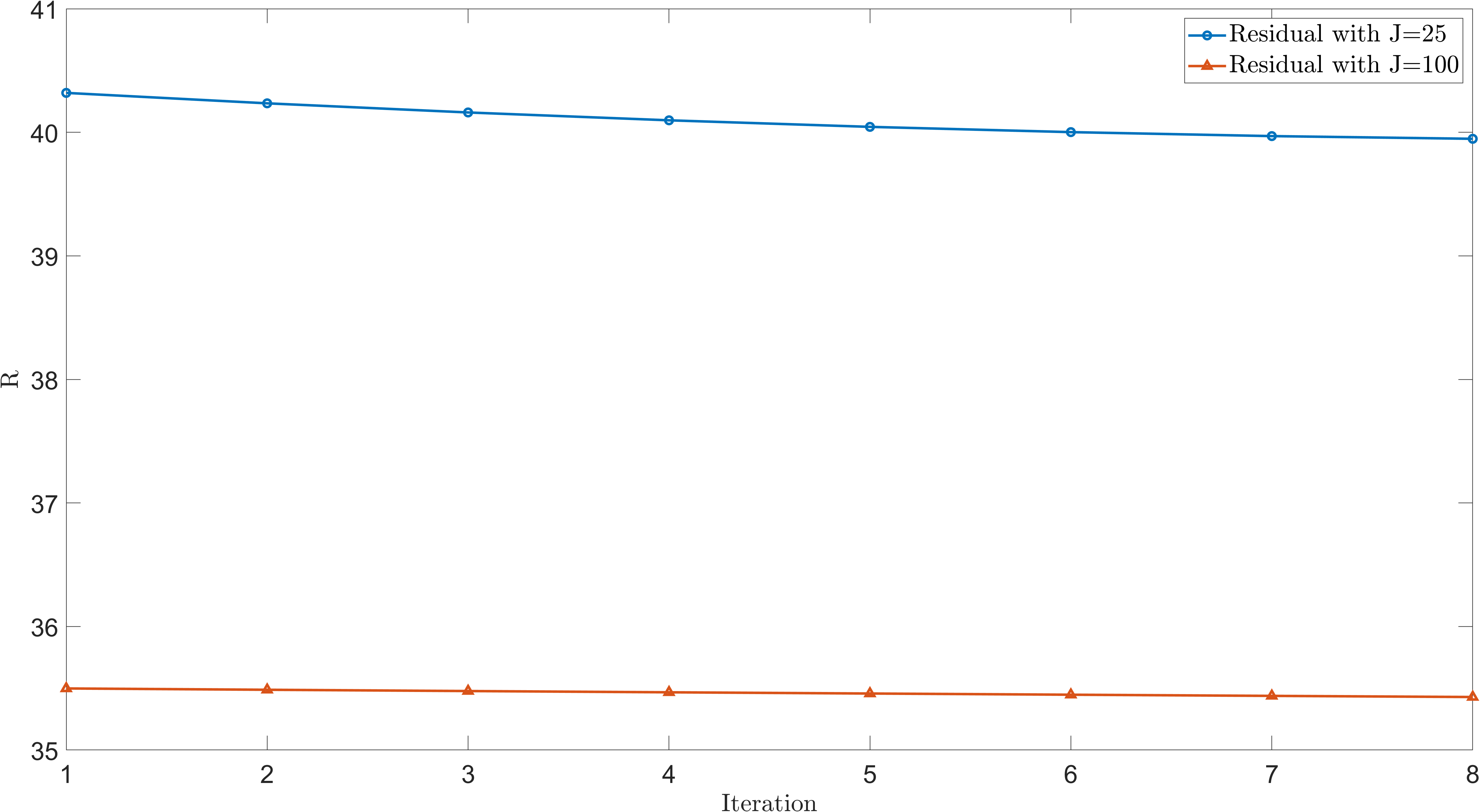}}\\
	\subfloat[]{\includegraphics[width=.49\textwidth]{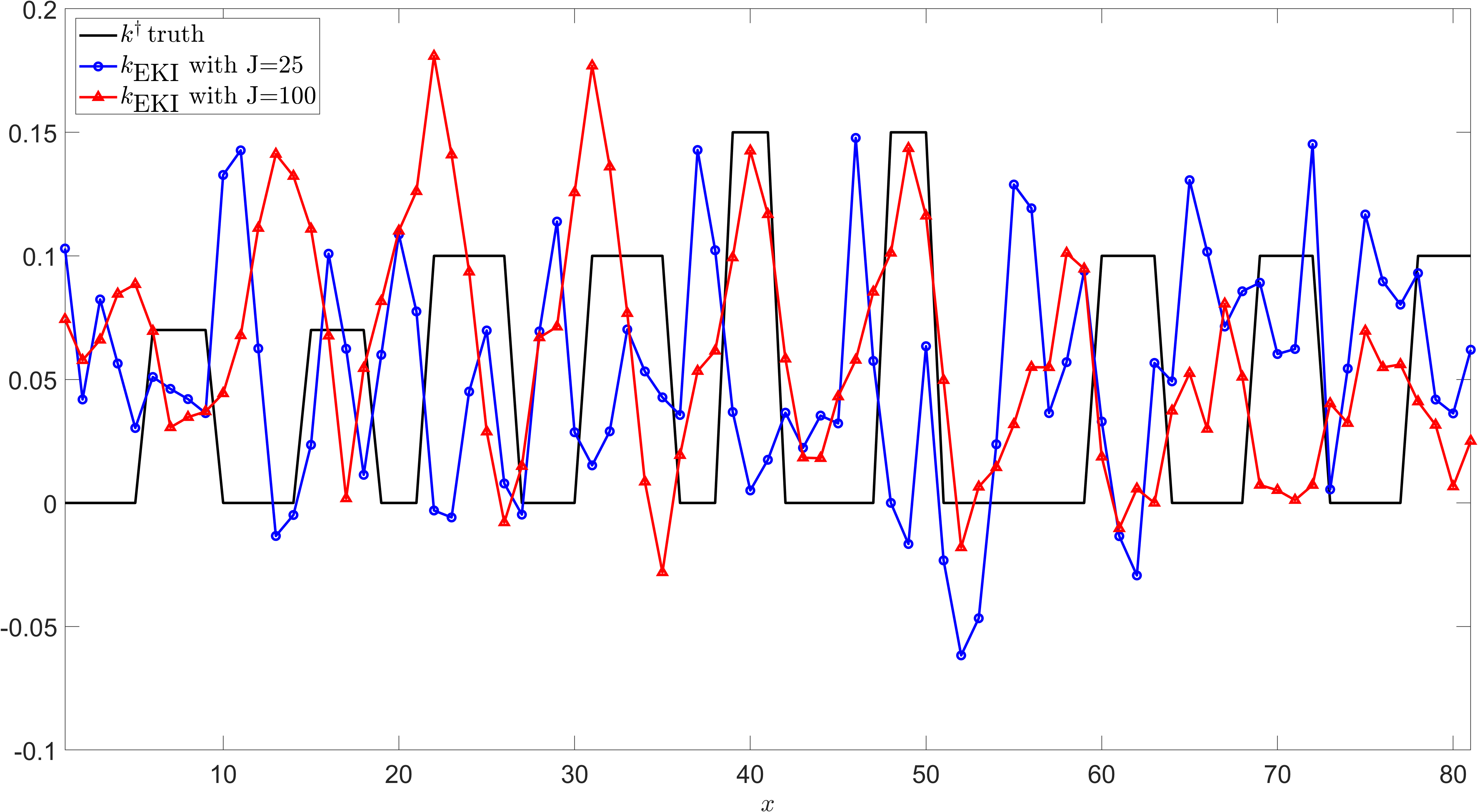}\hspace{0.1cm}}
	\subfloat[]{\includegraphics[width=.49\textwidth]{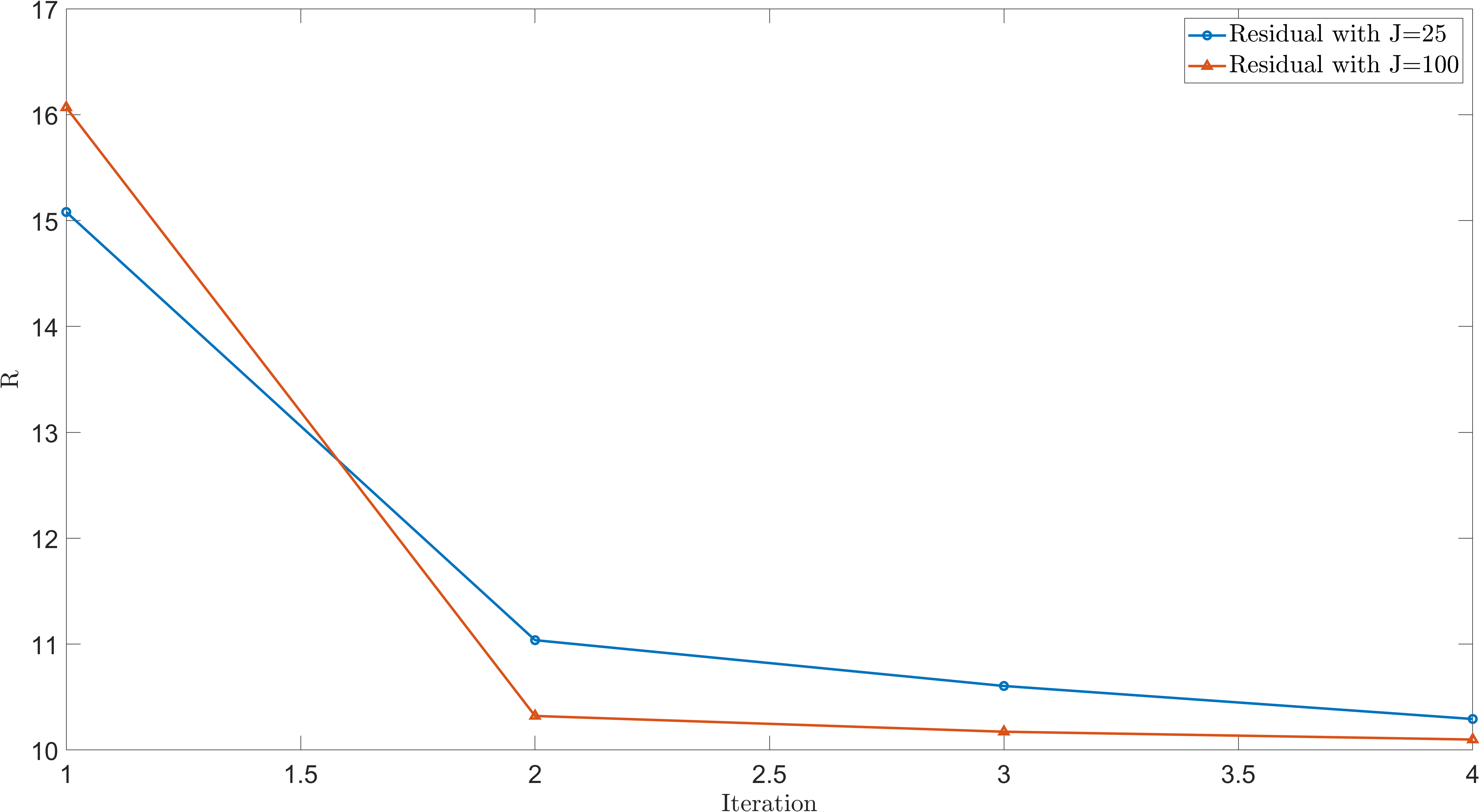}}\\
	\subfloat[]{\includegraphics[width=.49\textwidth]{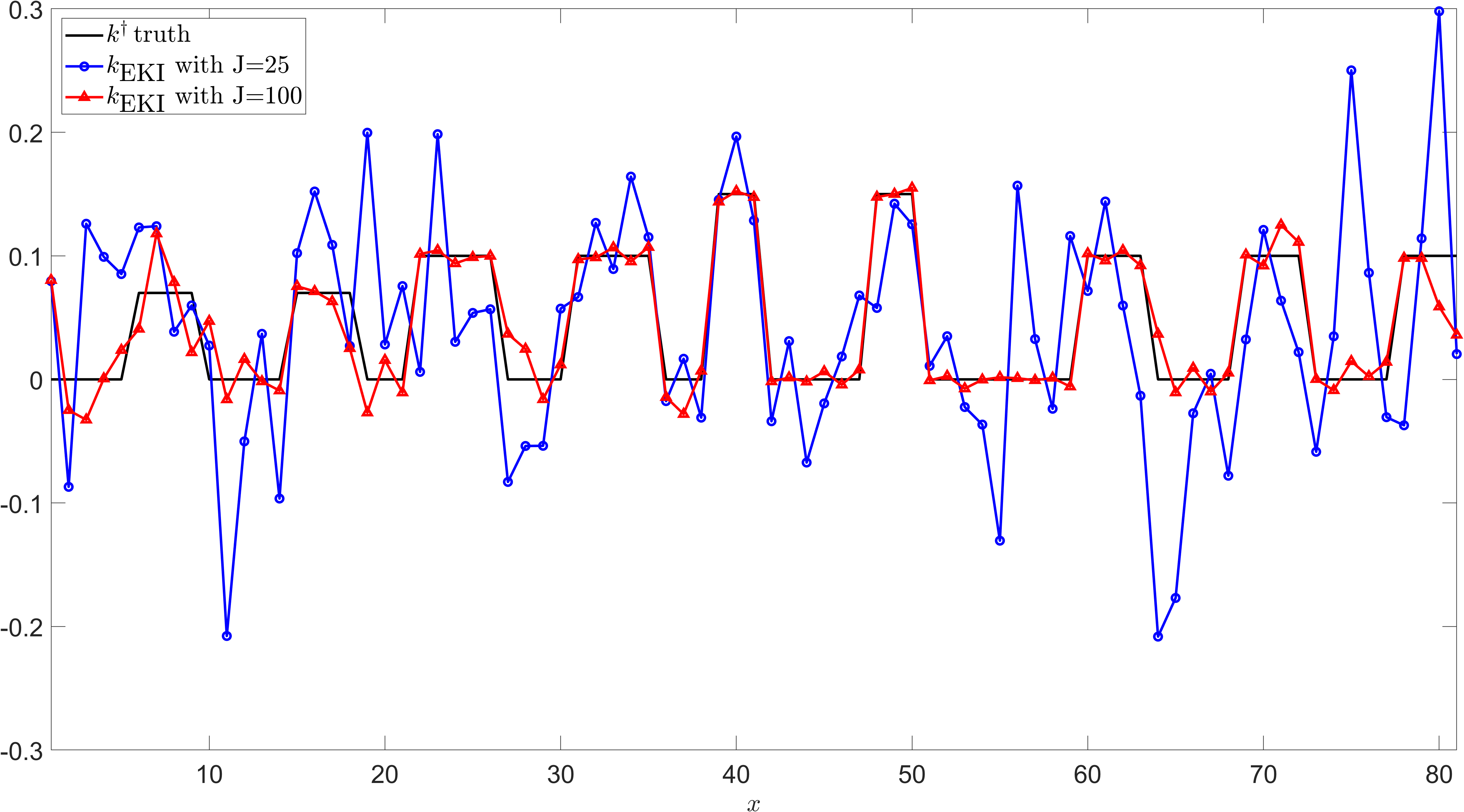}\hspace{0.1cm}}
	\subfloat[]{\includegraphics[width=.49\textwidth]{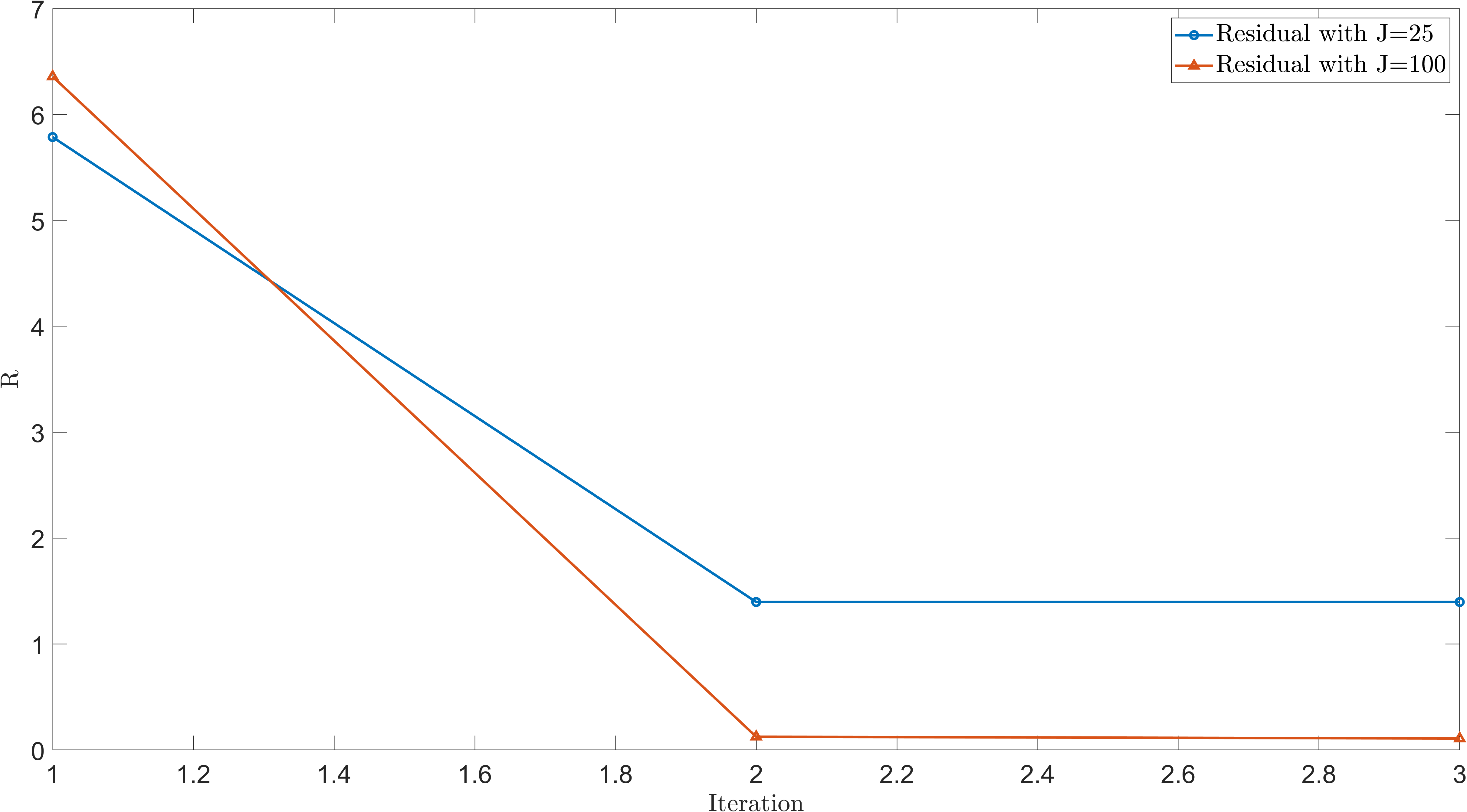}}\\
	\caption{\textit{Top line} (a),(b): with $\gamma=0.005$ and $\beta=6000$. \textit{Middle line} (c),(d): with $\gamma=10^{-7}$ and $\beta=3000$. \textit{Bottom line} (e),(f): with noise-free and $\beta=1000$. \textit{Left column} (a),(c),(e): comparison between the truth value of Winkler coefficient $k^\dagger$ and its reconstruction $k_{\textrm{EKI}}$ with $J=25,100$ particles. \textit{Right column} (b),(d),(f): residual of the Winkler coefficient $k$ for $J=25,1000$.}
	\label{residual_piecewise}
\end{figure}

\noindent\textbf{The inverse problem:}
Now let is analyze the inverse problem, that means, suppose we have the solution $w$ of the system (\ref{eq:biharmonic_sys}), we want to reconstruct the Winkler coefficient $k$, showing the residual that corresponds to the difference between the truth value $k^\dagger$ and the coefficient reconstructed by the EKI method $k_{\textrm{EKI}}$. We distinguish the case with noise and the case without noise and we observe that the EKI method in the absence of noise is easily able to reconstruct the Winkler coefficient $k$, while in the case in which the noise is considerable, i.e. $\gamma=0.005$, it finds some difficulty, but the residual, see the second column in fig. \ref{residual_piecewise}, decreases in both cases, so the EKI method is in any case a good estimator.
\begin{figure}[!ht]
	\centering
	\subfloat[]{
		\includegraphics[width=.32\textwidth]{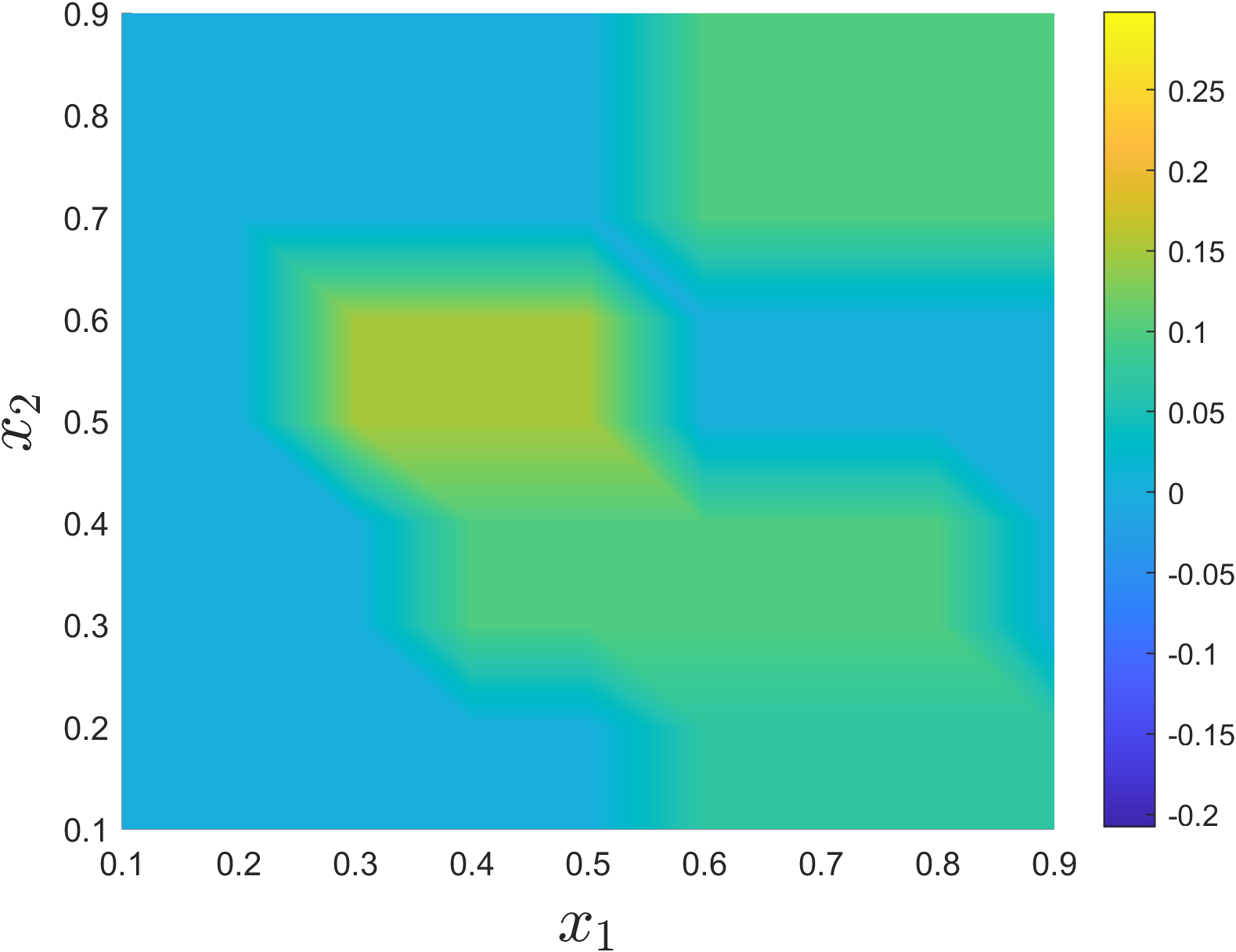}}\hspace{0.05cm}
	\subfloat[]{
		\includegraphics[width=.32\textwidth]{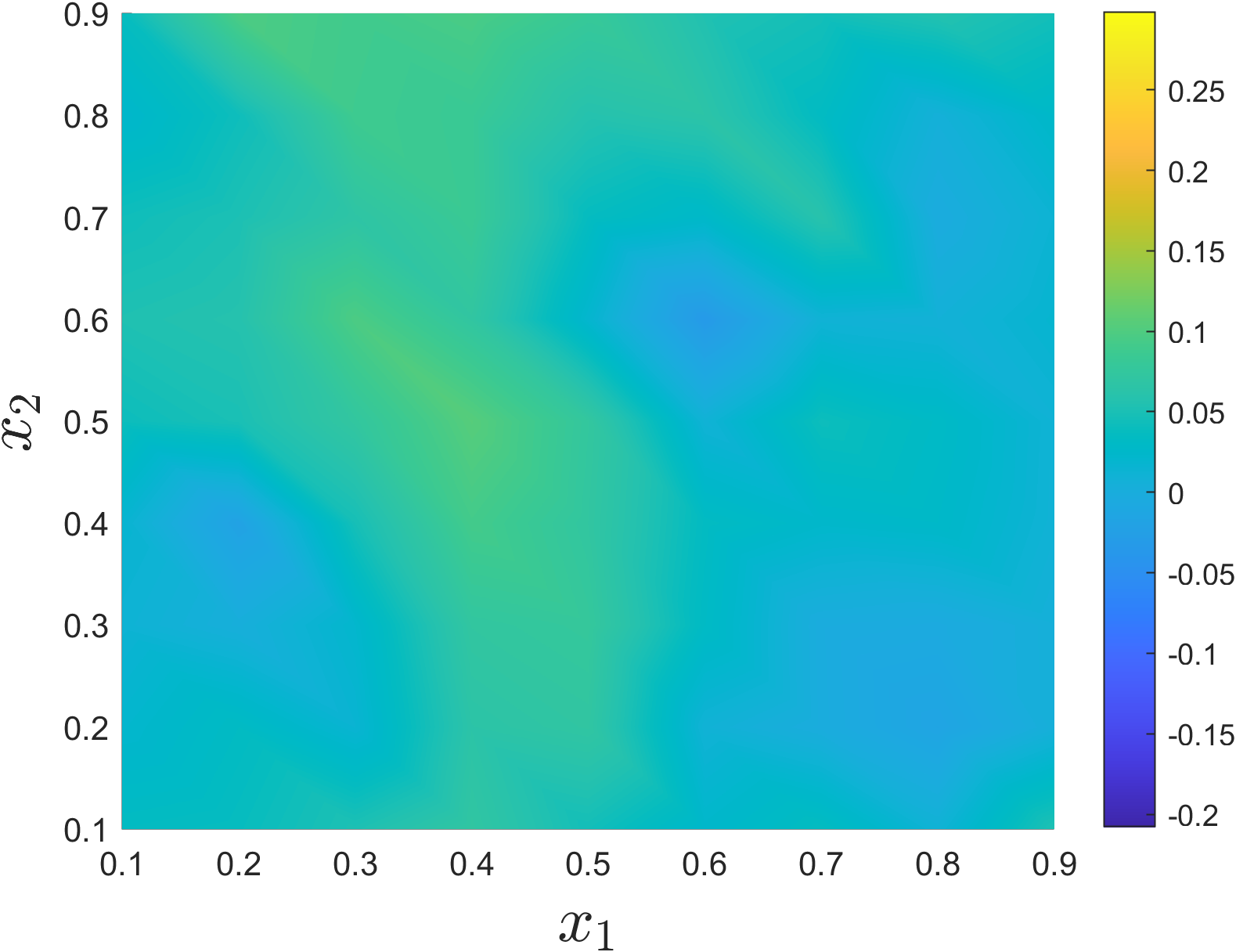}}\hspace{0.05cm}
	\subfloat[]{\includegraphics[width=.32\textwidth]{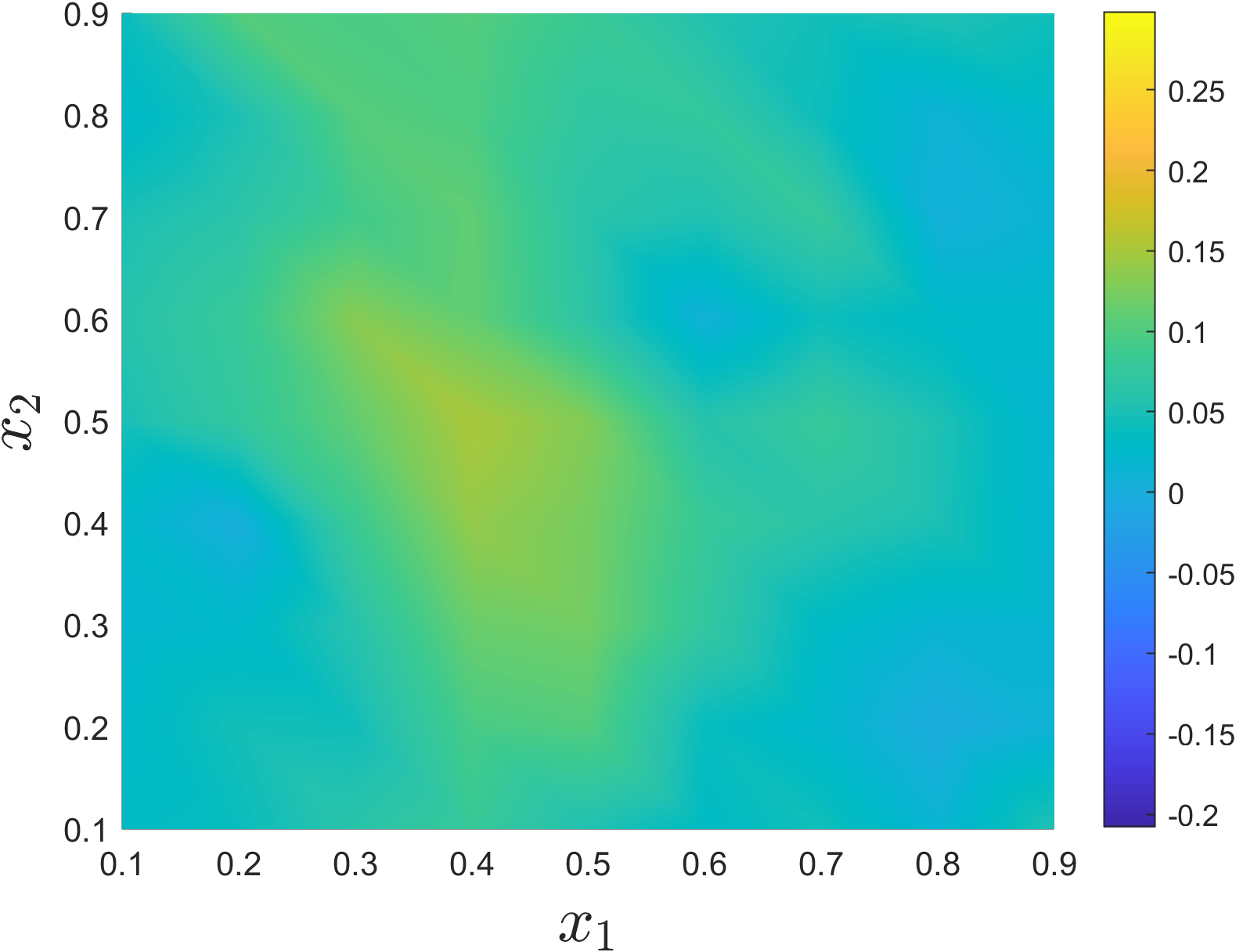}}\\
	\subfloat[]{
		\includegraphics[width=.32\textwidth]{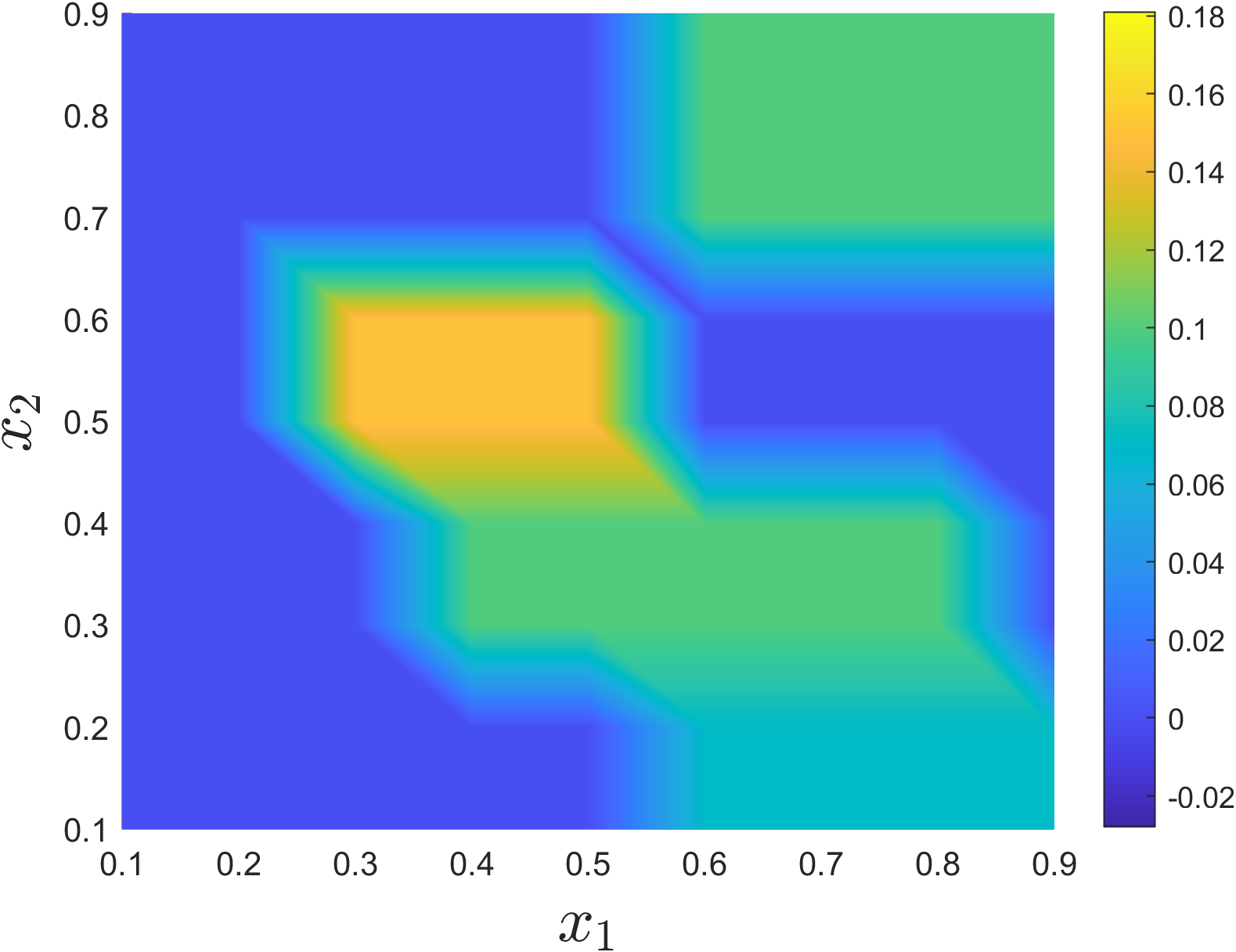}}\hspace{0.05cm}
	\subfloat[]{
		\includegraphics[width=.32\textwidth]{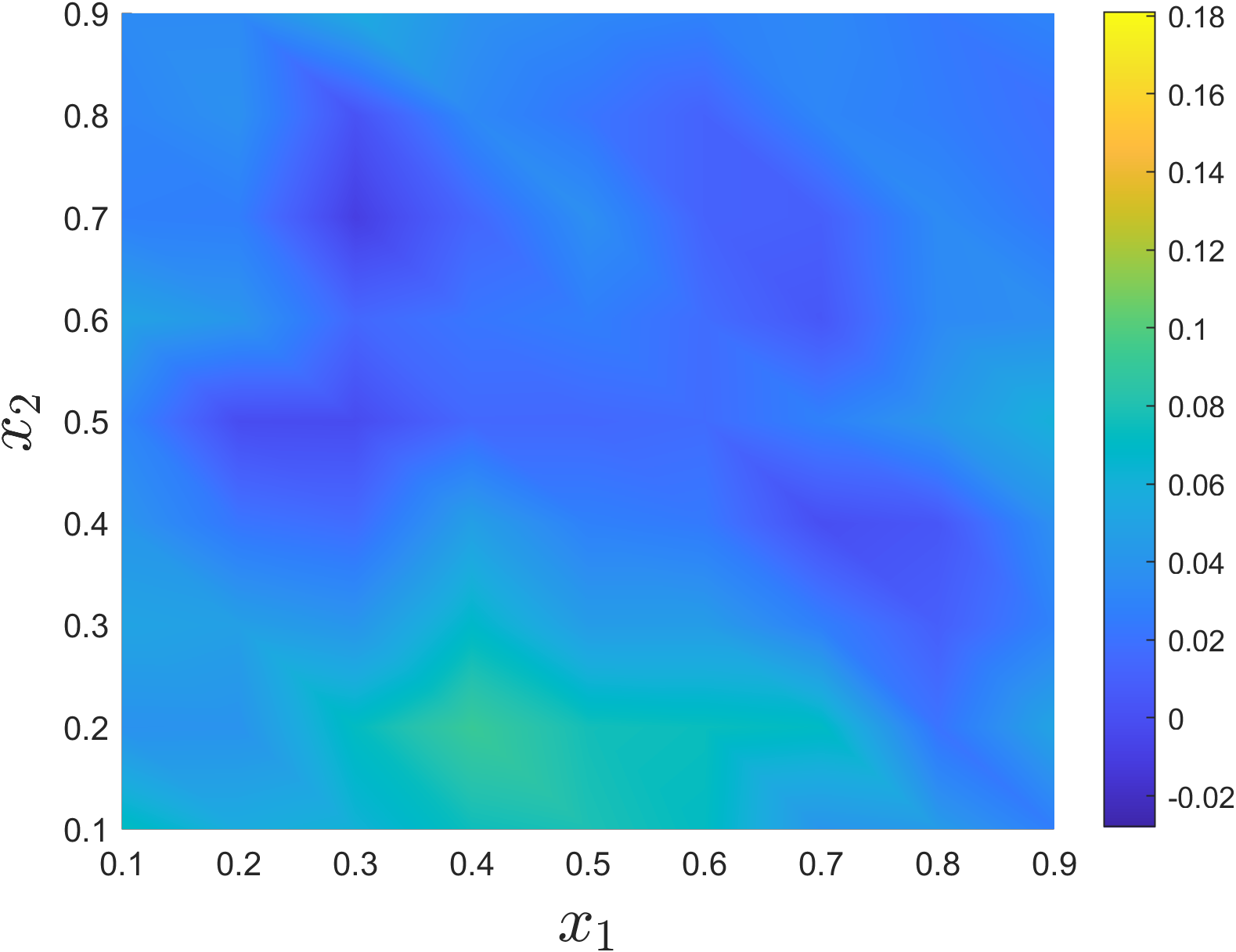}}\hspace{0.05cm}
	\subfloat[]{\includegraphics[width=.32\textwidth]{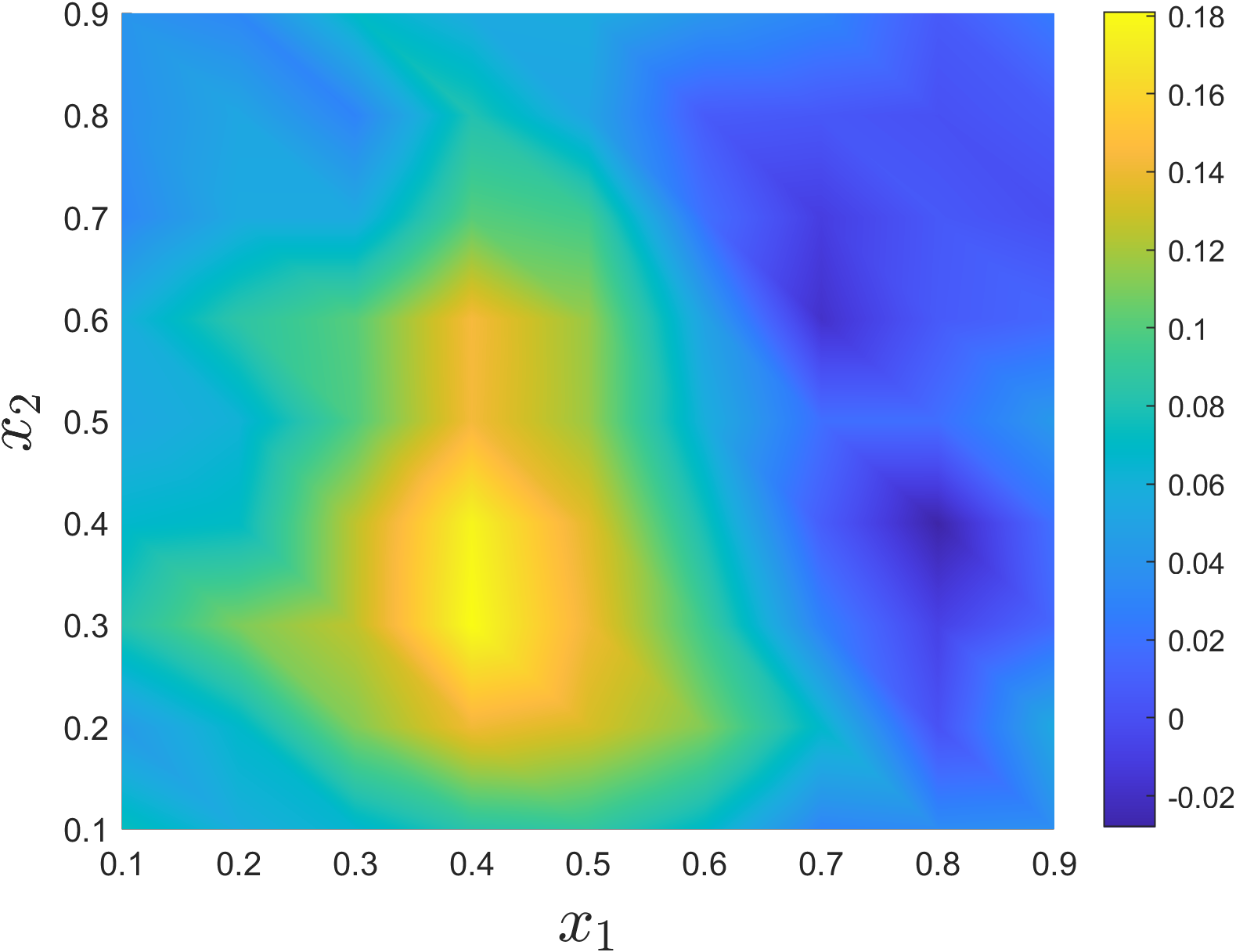}}\\
	\subfloat[]{
		\includegraphics[width=.32\textwidth]{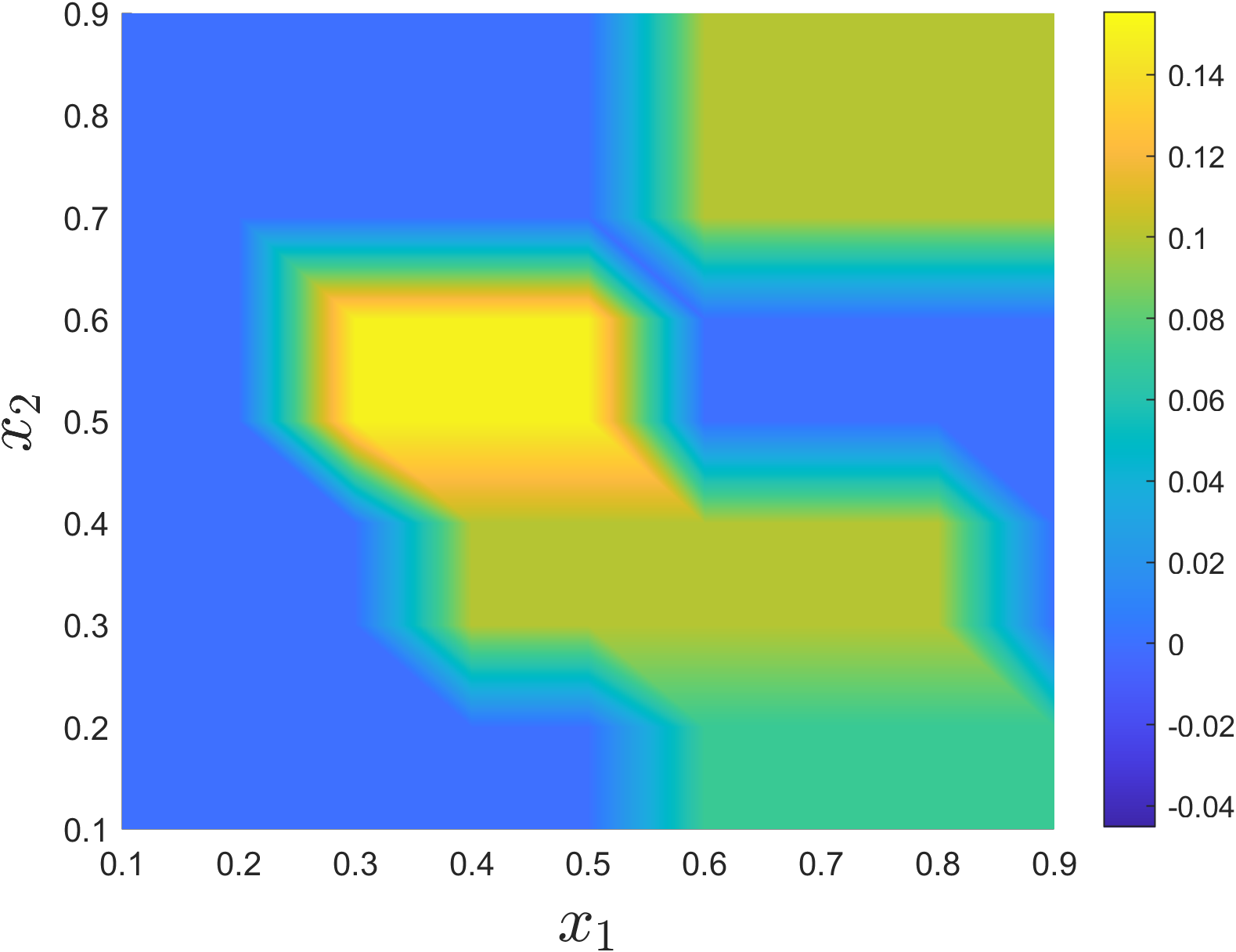}}\hspace{0.05cm}
	\subfloat[]{
		\includegraphics[width=.32\textwidth]{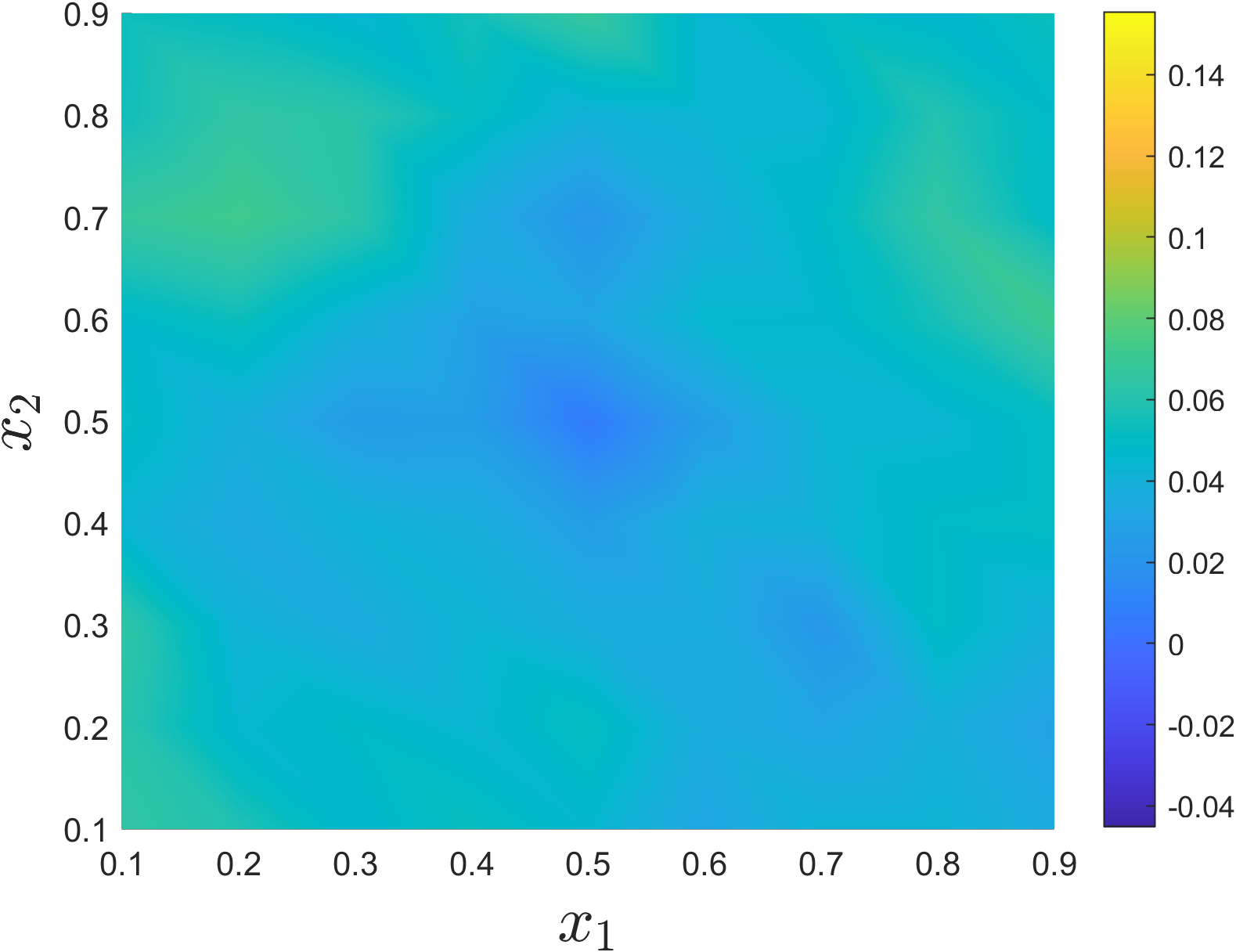}}\hspace{0.05cm}
	\subfloat[]{\includegraphics[width=.32\textwidth]{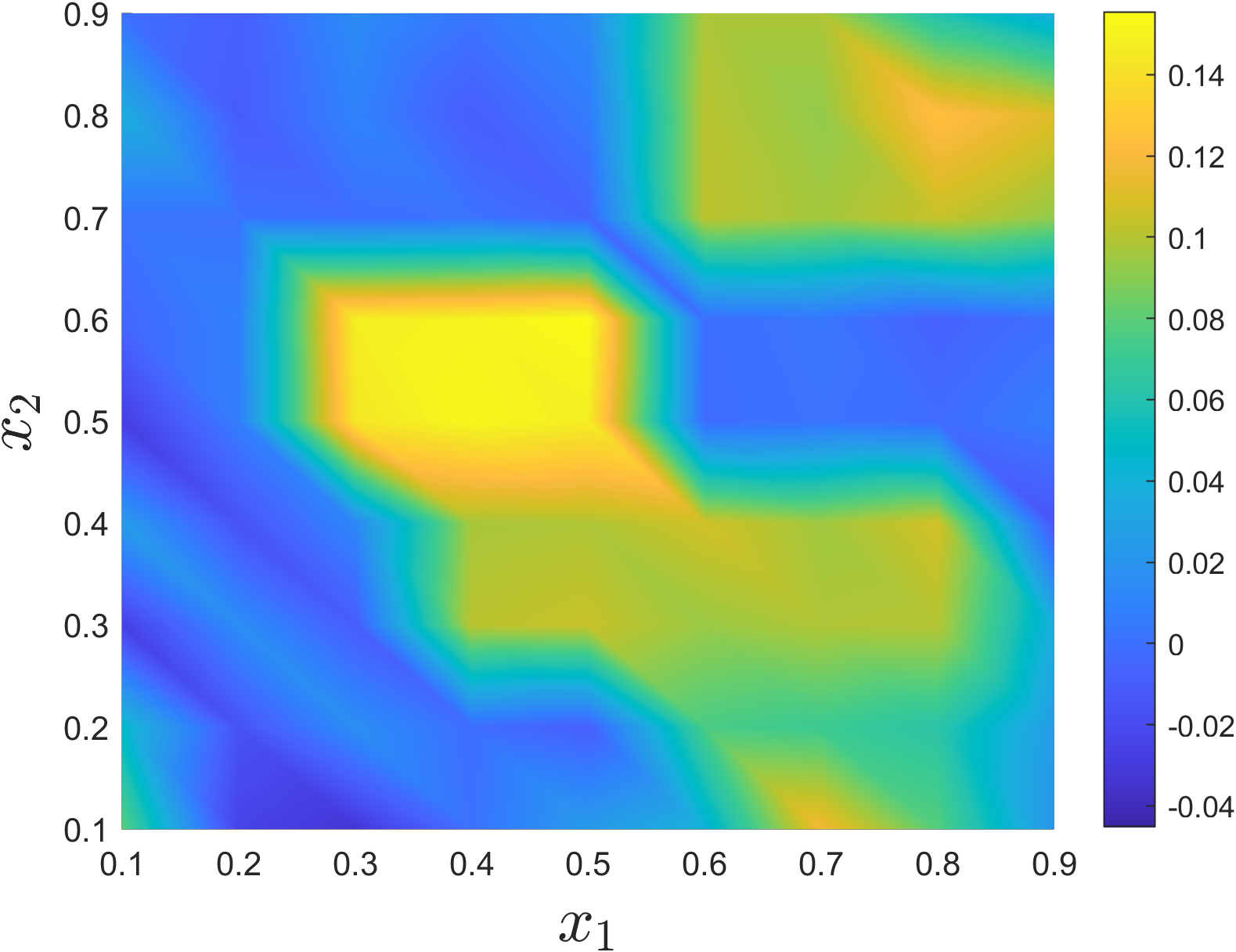}}\\
	\caption{We set $J=100$ number of particles. \textit{Top line} (a),(b),(c): with $\gamma=0.005$ and $\beta=6000$. \textit{Middle line} (d),(e),(f): with $\gamma=10^{-7}$ and $\beta=3000$. \textit{Bottom line} (g),(h),(i): with noise-free and $\beta=1000$. \textit{Left column} (a),(d),(g): Numerical solution of $k$. \textit{Middle column} (b),(e),(h): prior measure $\mu_0$ of $k$. \textit{Right column} (c),(f),(i): Reconstruction of $k$.}
\end{figure}
\noindent To better observe the behavior of the EKI algorithm applied to this case, we will try to plot, in the first column of the Fig. \ref{residual_piecewise}, in one-dimension.

\section{Summary and perspectives}\label{section_sum_perp}
In this paper we have introduced the EKI method to solve an inverse problem simplified by the Winkler model to be able to better reconstruct the subgrade reaction coefficient from the measurements of the transverse displacements induced by a load concentrated in a precise point of a thin plate installed in the foundation of an existing building. Considering that this is a ill-posed problem, we were able to reconstruct within our domain both in the case of the \textit{direct problem} (from k obtain w) and in the case of the \textit{inverse problem} (from w obtain k). The numerical results have shown that in the noise-free regime it is easier reconstruct the solution for the inverse problem with respect to the solution reconstructed in a noise regime, but in both cases the residual, i.e. the difference in $\Gamma$-norm between the truth value and the reconstructed value, goes to zero. As future perspective it would be interesting to use a denoising numerical method from the numerical solution and only then apply the optimal estimator EKI method.
\subsection{Acknowledgments}
Leonardo Scandurra is member of the ``National Group for Scientific Computation (GNCS-INDAM)'' and acknowledges support by MUR (Ministry of University and Research) PRIN2017 project number 201758MTR2-007 ``Direct and inverse problems for partial differential equations: theoretical aspects and applications''. Thanks to Edi Rosset and Eva Sincich, Associate Professors of Mathematical Analysis at the Università degli studi di Trieste for the helpful discussions and for providing me some useful material about the inverse problems.
\bibliographystyle{abbrv}
\bibliography{bibfile}

\end{document}